\font\steptwo=cmb10 scaled\magstep2
\font\stepthree=cmb10 scaled\magstep3

\magnification=\magstephalf
\settabs 18 \columns
\voffset-1cm
\hsize=16truecm

\def\b{\bigskip}
\def\bb{\bigskip\bigskip}

\def\no{\noindent}
\def\r{\rightline}
\def\ce{\centerline}
\def\ve{\vfill\eject}

\font\got=eufm10 scaled\magstep0

\def\g{{\got g}}
\def\sl{\hat\partial }
\def\cb{\hbox{${\cal B}$}}
\def\ba{${\cal B}$}
\def\p{\partial}
\def\ca{\hbox{${\cal E}$}}

\def\today{\ifcase\month\or January\or February\or
March\or April\or  May\or June\or July\or August\or
September\or October\or November\or  December\fi
\space\number\day, \number\year }

\r \today
\vskip1in

{\ce {\stepthree   THE IDEALS OF}}
\b
{\ce{\stepthree  FREE DIFFERENTIAL ALGEBRAS}}
\bb
{\ce {C. Fr\o nsdal}}
{\ce {Physics Department, University of California, Los
Angeles CA 90024, USA}}
\ce {and}
{\ce {A. Galindo}} {\ce {Departamento de F\'{\i}sica
Te\'orica, Universidad Complutense, 28040 Madrid, Spain}}
\bb

\no{\it ABSTRACT.}
We consider the free ${\bf C}$-algebra ${\cal B}_q$ with $N$
generators $\{\xi_i\}_{i = 1,...,N}$, together with a set of
$N$ differential operators $\{\p_i\}_{i = 1,...,N}$ that act
as twisted derivations on \ba$_q$ ~according to the rule
$\p_i\xi_j = \delta_{ij} + q_{ij}\xi_j\p_i$; that is, $
\forall x \in \cb_q, \partial_i(\xi_jx) =
\delta_{ij}x + q_{ij}\xi_j\partial_i x,$
and $\partial_i{\bf C} = 0$.  The suffix $q$ on $\cb_q$
stands for $\{q_{ij}\}_{i,j \in \{1,...,N\}}$ and is
interpreted as a point in parameter space, $q =
\{q_{ij}\}\in {\bf C}^{N^2}$.  A constant $C \in {\cal B}_q$
is a nontrivial element with the property $\p_iC = 0,~ i =
1,...,N$.  To each point in parameter space there correponds
a unique set of constants and a differential complex.  There
are no constants when the parameters $q_{ij}$ are in general
position.  We obtain some precise results concerning the
algebraic surfaces in parameter space on which constants
exist.  Let ${\cal I}_q$ denote the ideal generated by the
constants.  We relate the quotient algebras $\cb_q' =
\cb_q/{\cal I}_q$ to Yang-Baxter algebras and, in
particular, to quantized Kac-Moody algebras.  The
differential complex is a generalization of that of a
quantized Kac-Moody algebra described in terms of Serre
generators.  Integrability conditions for $q$-differential
equations are related to Hochschild cohomology.  It is shown
that $H^p(\cb_q',\cb_q') = 0$ for $p \geq 1$.  The intimate
relationship to generalized, quantized Kac-Moody algebras
suggests an approach to the problem of classification of
these algebras.

\ve
\no {\steptwo   Introduction.}
\b
\quad A recent study of the universal R-matrix of
quantum groups led to the study of a type of free
differential algebras. The positive Serre generators
generate an algebra with certain relations,
Drinfeld's  quantized Serre relations, and the
action of the other generators can be expressed in
terms of $q$-differentiation operators. But the
Serre relations can be replaced by others, and wide
generalizations are possible. This leads to the
study of free differential algebras and their
ideals. In this paper we obtain results that bear
on the classification of ideals in  quantum groups
and Kac-Moody algebras, and
on the classification of a new series of  quantum
Hopf algebras (generalized quantum groups).

To illustrate the type of generalization that is encompassed
here, consider the generalized Cartan matrices
$$
 \pmatrix{2&-1&n\cr-1&2&-1\cr n&-1&2\cr},\quad n = -1,0,2.
$$
For $n = 0$ it corresponds to a simple Lie algebra, and to a
finite quantum group.  For $n = -1$ it is the Cartan matrix
of a Kac-Moody algebra of affine type, also quantizable.
The quantized Serre ideals include generators of the type
$$
\eqalign{
n = 0&:  \quad [e_1,e_3]_q:=e_1e_3-qe_3e_1,
\cr
n = -1&: \quad [e_1,[e_1,e_3]_q]_{q'}.
}
$$
For $ n=2$, there is no generator of this type, constructed
from $e_1$ and $e_3$, but instead there is a generator of a
new type, namely
$$
n = 2:\quad [[e_2,e_1]_q,e_3]_{q'} + [[e_2,e_3]_q,e_1]_{q'}.
$$
The parameters $q,q'$ tend to 1 in the ``classical" limit.
As this example shows, generalized quantized Kac-Moody
algebras of indefinite type are characterized by unexpected
Serre-type ideals.

Section 1 begins with the abstract definition of a family of
free algebras with $q$-differential structure.  The
connection with quantum groups and universal R-matrices is
reviewed in Section 2.  We demonstrate the advantage of our
methods by evaluating the highest root vectors for
$U_q(A_l)$ (already found by Jimbo [J]) and for $U_q(C_l)$.

Let $\cb$ be the ${\bf C}$-algebra freely generated by
$\xi_1,\ldots,\xi_N$.  A $q$-differential structure is a set
of operators $\partial_i,\ldots,\partial_N$ that act on
$\cb$ by the rule
$$
\partial_i(\xi_jx) = \partial_{ij}x + q_{ij}\xi_j\partial_ix,
\quad \partial_i{\bf C} = 0.
$$
This action involves a set $\{q_{ij}\}$ of complex
parameters.  We denote by $\cb_q$ the algebra $\cb$ endowed
with the $q$-differential structure.  An interesting
property of $\cb_q$ is that, when the parameters are in
general position, $q$-differential equations of the form
$\partial_ix=y_i$, $i=1,\ldots,N$, can be solved for $x$, for
any choice of the ``one-form" $\{y_i\}$; all one-forms are
exact.  Let us call ``exceptional" those points in parameter
space for which this is not true.  These are precisely the
same as those for which there exist homogeneous elements in
$\cb_q$, of degree higher than zero, satisfying $\partial_iC
= 0$, $i=1,\ldots,N$; such elements are called {\it
constants}.  The constants in $\cb_q$ generate an ideal
${\cal I}_q$ in $\cb_q$ and allows to define the quotient
algebra $\cb_q^\prime = \cb_q/{\cal I}_q$.  Quantized
Kac-Moody algebras are a particular case that is the subject
of Section 2.  A main goal is the complete classification of
all special points in parameter space; in the following
sense.  Two points $q$ and $q'$ in the space of parameters
are said to be equivalent if the ideals ${\cal I}_q$ and
${\cal I}_{q'}$ of $\cb_q$ and $\cb_{q'}$ coincide (as
subalgebras of $\cb$).  An alternative and probably more
fundamental classification, where equivalence is instead
based on isomorphism of the respective quotient algebras, is
not contemplated in this paper.  The complete classification
of $\cb'_q$-algebras would provide, in particular, a partial
classification of quantized Kac-Moody algebras of the most
general type.
\b

The existence of constants in $\cb_q$ is revealed by the
reduction in rank, at exceptional points in parameter space,
of the matrix $S$ defined by
$$
 S^{j_1...j_p}_{i_1...i_p}
:=\partial_{i_p}\ldots\partial_{i_1}(\xi_{j_1}...\xi_{j_p}).
$$
\no The projection of this matrix on $\cb_q^\prime$ is
invertible and the components of the inverse matrix appear
as coefficients of the expansion of the universal R-matrix
in terms of Serre generators.  Here we make contact with the
work of Varchenko [V] on quantum groups.

An action of $\cb_q'$ on $\cb_q'$ is defined via the
homomorphism that sends $\xi_i$ to $\partial_i$, $i =
1,...,N$.  It is shown that $H^p(\cb_q',\cb_q')$, defined
via this action, is zero for $p \geq 1$.  The meaning of
this result in the case of quantized Kac-Moody algebras is
as follows.  Let $\hat{\cal A}_+ $ be the algebra freely
generated by $\{e_i\}_{i = 1,...,N}$, identified with
$\cb_q$.  Serre type ideals in $\hat{\cal A}_+ $ are
generated by elements annihilated by the $f_i$'s; they are
precisely the constants in $\cb_q$.  One can therefore draw
the conclusion that the generalized, quantized Kac-Moody
algebra associated to $\cb'_q$ is rigid with respect to
deformations that respect the Cartan decomposition.

Triviality of cohomology in this algebraic setting is not
surprising.  Interesting, nontrivial cohomology depends on
completion.  A particular type of completion is implied by
the interpretation of the differential structure in terms of
finite difference operators, possibly related to difference
equations of the type studied by Smirnov [S][FR].

The advantages of the Serre formulation of (quantized)
Kac-Moody algebras are obvious.  It is natural to ask
whether some aspects of the cohomology of these algebras (as
generalized) can be formulated in terms of linear forms on
the span of the Serre generators.  We show, in Section 3.5,
that the natural definition of ``Serre cohomology" is in
terms of multilinear $p$-cochains restricted to closed
chains.  For cochains valued in $\cb'_q$ the differential is
the projection on closed chains of the map that sends a
$p$-cochain $z$ to the $(p+1)$-cochain $dz$ according to the
formula
$$
dz(\xi_{i_0}\otimes \xi_{i_1}\otimes...\otimes\xi_{i_p}) =
\partial_{i_0}z(\xi_{i_1}\otimes...\otimes\xi_{i_p}).
$$
The cohomology is not trivial, but $z$ is exact iff it is
``strongly closed"; that is, iff\footnote*{Here, and in the
formulas that follow, the summation is over all repeated
indices, running independently over $\{1,...,N\}$.  We use a
multi-index notation where ${\underline i}$ stands for
$i_1,...i_k$.}
$$
\sum_{{\underline i}}C^{{\underline i}}\,\partial_{i_1}...\partial_{i_k} =
0\quad \Rightarrow \quad
\sum_{{\underline i}}C^{{\underline
i}}\,\partial_{i_1}...\partial_{i_{k-1}}z(\xi_{i_k}\otimes
\xi_{j_1}\otimes...\otimes \xi_{j_p}) = 0.
$$

We must explain this statement.  If there are coefficients
$C^{{\underline i}}$ such that
$$
\sum_{{\underline i}}C^{{\underline i}} \,\xi_{i_1}\ldots\xi_{i_k}
:= \sum_{i_1,...,i_k}C^{i_1\ldots
i_k} \xi_{i_1}\ldots\xi_{i_k} \in \cb_q
$$
is a constant, then
$$
\sum_{{\underline i}} C^{{\underline i}}\, \xi_{i_1}\ldots\xi_{i_k} = 0
$$
is a relation in $\cb_q^\prime$, and (as is shown in
subsection 1.3.3) the operator
$$
\sum_{{\underline i}}
C^{{\underline i}}\,\partial_{i_1}\ldots\partial_{i_k}
$$
is identically zero.  The cohomology is nontrivial if there
are irreducible constants of polynomial order higher than 2.
For more details please turn to Section 3.5.

Section 4 contains results pertaining to the classification
of the algebras $\cb_q^\prime$; more precisely, we try to
determine the exceptional points in parameter space at which
constants occur in $\cb_q$.  Let $Q = (1,\ldots,n)$ and $Q_i
= (1\ldots\hat{\hbox{\it\i}}\ldots n)$.  Let $\cb_{1...n} =
\cb_Q$ be the space of polynomials linear in
$\xi_1,\ldots,\xi_n$ separately.  It is shown that, if the
parameters associated with $Q_i$, $i=1,\ldots,n$, namely
$\{q_{jk},~\{j,k\}\subset Q_i\}$, are in general position;
more precisely, if there are no constants in $\cb_{Q_i}, i =
1,...,n$, then constants exist in $\cb_{1...n}$ if and only
if
$$
1 - \prod_{i\not=j\in Q} q_{ij}  = 0~,
\eqno(1)
$$
\no and that the dimension of the subspace of constants in
$\cb_{1...n}$ is then $(n-2)!$.  Analogous results, for
arbitrary sets $Q$ (with repetitions), will be reported
elsewhere.

Physical applications of hyperbolic Kac-Moody algebras
appear in connection with dimensional reduction of general
relativity [N].  In other contexts it is interesting to look
for finite-dimensional representations.  It seems likely
that finite representations of quantized Kac-Moody algebras
of hyperbolic or more general type exist only for parameters
at roots of unity.  It is interesting to notice that all the
constraints turn out to imply a factorization of unity.  In
this connection it may be productive to take another point
of view.  Instead of regarding the $q_{ij}$ as complex
parameters, one may regard them as generators of a
commutative algebra, and replace the field ${\bf C}$ by the
ring of polynomials in $q_{ij},i,j = 1,...,n$.  In this
interpretation there is a unique algebra generated by
$q_{ij}, \xi_i$, and the left hand side of Eq.(1) generates
an ideal.  The problem is then one of classification of
certain ideals of a commutative algebra.  Compare [V].

Completion of the work contained in Section 4 would go some
way towards the classification of the algebras
$\cb_q^\prime$.  For more complete results we expect that
geometrical methods, such as those of Varchenko [V], may be
the most powerful.  We suggest, in particular, that a study
of the holonomy of the arrangement of surfaces defined by
the exceptional points in parameter space may be useful and
interesting.

Section 5 gives a complete account of constants in the
subspace $\cb_{123}$ of polynomials separately linear in
three generators.

\bb
\no
{\steptwo 1. q-differential algebras.}
\b

In this section we present the principal players: the freely
generated algebra {\cb} with its $q$-differential structure
and the symmetric form $S$, with their most basic
properties.

\b
\no
{\bf 1.1.  Free algebra and $q$-differential structure.}

\b \no{\it1.1.1.} Let {\cb} denote the unital {\bf C}-algebra
freely generated by $\xi_1,...,\xi_N$, with its natural
grading, $\cb = \bigoplus_{n\geq 0} \cb(n)$, where $\cb(n)$
contains homogeneous polynomials of degree $n$.  Suppose
given a map
$$
q: [N] \times [N] \rightarrow {\bf C},\quad (i,j) \mapsto q_{ij},
\quad [N]:=\{1,...,N\}.
$$
The $q_{ij}$'s are {\it the parameters}; a choice of
parameters will be interpreted as selecting a point in the
space $V = {\bf C}^{N^2}$.

\b
\no
{\it 1.1.2.} For a fixed
choice of parameters, let $\partial_1,...\partial_N$ be the
set of linear, $q$-differential operators
$$
\partial_i:  \cb(n) \rightarrow \cb(n-1),~ n\geq 1,
\quad \cb(0) \rightarrow
0,\quad  i = 1,...,N,
$$
defined by
$$
\partial_i(\xi_jx) = \delta_{ij}x + q_{ij}\xi_j\partial_ix,
\quad \forall  x\in \cb.
$$
In particular,
$$
\partial_i(\xi^r_i) = [r]_{q_{ii}}\,\xi_i^{r-1},\quad
[r]_{q} := 1 + q + ...
+q^{r-1}.
$$
\no Let $\cb_q$ denote the algebra {\cb} endowed with this
$q$-differential structure.  Thus $\cb_q$, as an algebra, is
identified with $\cb$.

\b
\no
{\it 1.1.3.} Let $\widehat{\cb}$ denote the unital {\bf
C}-algebra freely generated by $\sl_1,...,\sl_N$, with the
natural grading.  There is a unique homomorphism
$$
D_q: \widehat{\cb} \rightarrow {\rm End }\,\cb,
$$
such that $ \sl_i \mapsto \partial_i, i = 1,...,N$.  Let
$\cb^*_q$ denote the image of $\widehat{\cb}$ in End$\,$\cb;
it is the algebra of linear differential operators generated
by $\partial_1,...,\partial_N$, but unlike $\widehat{\cb}$
it is not always freely generated.

\b
\no
{\it 1.1.4.  Example.} If $q_{12}q_{21} = 1$, then
$(\partial_1\partial_2 - q_{21}\partial_2\partial_1)x = 0$,
for all $x$ in $\cb_q$, and hence $\partial_1\partial_2 -
q_{21}\partial_2\partial_1 = 0$ is a relation of $\cb^*_q$.
{\it Proof.} For all $y \in \cb_q$, if $ q_{12}q_{21} = 1$,
we have
$$
\eqalign{
&(\partial_1\partial_2 -
q_{21}\partial_2\partial_1)(\xi_iy) =
q_{1i}q_{2i}\xi_i(\partial_1\partial_2 -
q_{21}\partial_2\partial_1)y,~ i = 1,2,...,N.
}
$$
Hence for all homogeneous $x\in \cb_q$ there is $f_x\in {\bf
C}$ such that
$$
(\partial_1\partial_2 - q_{21}\partial_2\partial_1)x =
f_xx(\partial_1\partial_2 - q_{21}\partial_2\partial_1)1 = 0.
$$

\b \no{\it 1.1.5.  Example.} If $q_{11} \neq 1$ is an $n$'th
root of unity, then $\partial_1^nx = 0$ for all $x$ in
$\cb_q$.

\b \no{\bf 1.1.6.  Theorem.} For generic $q$, $D_q$ is an
isomorphism.  {\it The exceptional values of} $q$, for which
Ker$D_q \neq \{0\}$, are determined by polynomial equations
in $q_{ii}, i = 1,...,N$ and in $\sigma_{ij} = q_{ij}q_{ji},
i \leq j = 1,...,N$.

The proof is in 1.3.

\b
\no{\bf 1.2. Constants.}

\b
To prepare for the proof of the theorem we introduce the
constants of $\cb_q$.

\b \no{\it 1.2.1.  Definition.} For a
fixed set of parameters, let $\cb_q = \bigoplus_{n\geq
0}\cb_q(n)$ denote the algebra \cb~ endowed with the
differential structure 1.1.2.  A {\it constant} in $\cb_q$
is an element $C \neq 0$ in $\bigoplus_{n \geq 1} \cb_q(n)$
that satisfies $\partial_iC = 0, i = 1,...,N$.  A {\it
homogeneous constant} is a constant that is homogeneous in
each generator.

Note the exclusion of $\cb(0)$.  Every constant is a sum of
homogeneous constants. The simplest examples are (1) if $q_{11} = -1$, then
$\xi_1^2$ is a constant, (2)
if $q_{12}q_{21} = 1$, then $\xi_1\xi_2 - q_{21}\xi_2\xi_1$ is a constant.

\b \no{\it 1.2.2.  Lemma.} For $q$ in general position there
are no constants in $\cb_q$.  The set of exceptional points
in $V$ is determined by polynomial equations in the $q_{ij},
i,j = 1,...,N$.

\b \no{\it Proof.} It is enough to consider homogeneous
constants.  A homogeneous constant is an element
$$
C = \sum_{{\underline i}} X^{{\underline
i}}\,\xi_{i_1}...\xi_{i_n} \in \cb_q(n),
$$
with complex coefficients $X^{{\underline i}}$, where the
sum is over effective permutations of the indices.  The
condition $\partial_iC=0, i = 1,...,N$ is equivalent to
$$
MC = 0,\quad M := \sum_{i = 1}^N\xi_i\partial_i.
$$
On the monomial basis for $\cb_q(n)$, the operator $M$ is a
matrix with coefficients that are monomials in the
$q_{ij}$'s.  Nonzero solutions exist if and only if the
determinant of this matrix is equal to zero.  This condition
is a polynomial equation in the $q_{ij}$'s.

\b \no{\it 1.2.3.  Lemma.} The exceptional set in 1.2.2 is
determined by polynomial equations in $q_{ii}, i = 1,...,N$,
and in $\sigma_{ij} = q_{ij}q_{ji}, i \leq j = 1,...,N$.

\b
\no{\it Proof.} A change of monomial basis leads to
$$
C = \sum_{\underline i} Y^{{\underline i}}\;{\xi_{i_1}...\xi_{i_n}\over
\prod q_{i_ki_l}},
$$
where the product includes a factor $q_{ij}$ for each
occurrance of the pair $(i,j),i<j$, in the order $j$ before
$i$, in the index set ${\underline i} = i_1,...,i_n$.  The
conditions $\partial_iC = 0$ reduce to linear equations for
the coefficients $Y^{{\underline i}}$.  After multiplying
each equation by a monomial in the $q_{ij}$'s, one obtains a
set of equations that involve only the $q_{ii}$ and the
$\sigma_{ij}$'s.  \b

\no
{\bf 1.3. Proof of the theorem 1.1.6.}

\b \no{\it 1.3.1.} Fix the parameters $q = \{q_{ij}\}$ and
let $\cb_q$ denote the algebra $\cb$ endowed with the
differential structure 1.1.2.  Let $\partial^{\rm t}_i, i =
1,...,N$, be another set of differential operators, defined
in the same way but with $q$ replaced by $q^{\rm t}$, where
$q^{\rm t}_{ij}:=q_{ji}$ for $i,j = 1,...,N$.  Let
$\cb_{q^{\rm t}}$ denote the algebra \cb~ endowed with this
new differential structure.

\b \no{\it 1.3.2.} Let $\widehat{\cb}\cb_q$ be the
universal, unital algebra generated by
$\xi_1,...,\xi_N,\sl_1,...,\sl_N$, with 
relations\footnote*{This algebra appears in work of Lusztig [L] and
Kashiwara [K], in the context of quantized Kac-Moody
algebras.}
$$
\sl_i\xi_j = \delta_{ij} + q_{ij}\xi_j\sl_i,
$$
Let $\widehat{\cb}\cb_{q^{\rm t}}$ be the same, but with
$q_{ij}$ replaced by $q_{ji}$.  There is a unique
anti-isomorphism
$$
\Phi:  \widehat{\cb}\cb_q \rightarrow
\widehat{\cb}\cb_{q^{\rm t}},
$$
such that $\xi_i \mapsto \sl^{\rm t}_i, \sl_i \mapsto \xi_i,
i = 1,...,N$.

\b \no{\it 1.3.3.} Now suppose that there is a homogeneous
constant in $\cb_q$.  Then $q = \{q_{ij}\}$ is exceptional,
and so is $q^{\rm t} = \{q^t_{ij}\}, q^t_{ij} = q_{ji}$, by
1.2.3.  Hence there is a homogeneous constant $C \in
\cb_{q^{\rm t}}$.  This implies that there are $f_i \in {\bf
C}, i = 1,...,N$, such that
$$
\sl_i^{\rm t}C = f_iC\sl_i^{\rm t} \in
\widehat{\cb}\cb_{q^{\rm t}}.
$$
Applying $\Phi^{-1}$ one gets
$$
\Phi^{-1}(C)\xi_i = f_i\xi_i\Phi^{-1}(C) \in
\widehat{\cb}\cb_q.
$$
This implies that, $\forall x \in \cb_q$, homogeneous in
each variable, there is $f_x \in {\bf C}$ such that
$$
\Phi^{-1}(C)x = f_xx\Phi^{-1}(C) \in  \widehat{\cb}\cb_{q};
$$
hence $\Phi^{-1}(C) \neq 0$ belongs to Ker$D_q$.
Conversely, any element of Ker$D_q$ is a sum of homogeneous
elements.  Let $\hat C'$ be a homogeneous element of
Ker$D_q$, of total degree $p$.  There are complex
coefficients $f'_i, i = 1,...,N$, such that $\hat C'\xi_i =
\hat C'_i + f'_i\xi_i\hat C' \in \widehat{\cb}\cb_{q}$, with
$\hat C'_i \in \,$Ker$\,D_q$ of total order $p-1$.
Iterating this by evaluating $\hat C'_i\xi_j$ and so on, one
eventually obtains an element $\hat C'\in\, $Ker$\,D_q$, of
order $m \geq 2$, (since there are no elements of order 1 in
Ker$D_q$), such that $\hat C'\xi_i = f_i\xi_i\hat C' \in
\widehat{\cb}\cb_q$.

Applying $\Phi$ one obtains $\sl_i^{\rm t}\Phi(\hat C') =
f_i\Phi(\hat C')\sl^{\rm t}_i$ and thus $\partial^{\rm
t}_i\Phi(\hat C') = 0$.  By 1.2.3 it follows that there is $C
\neq 0$ in $\cb_q(m)$ such that $\partial_iC = 0, i =
1,...,N$.  Hence $q$ is exceptional, in the sense of both
lemmas, and the theorem is proved.

\b \no{\it 1.3.4.  Comments.} The proof of the direct part
of the theorem demonstrates the existence of an monomorphism
from the space of constants in $\cb_q$ into Ker$D_q$.  In
the proof of the converse, however, the reduction to
elements of minimal degree is essential; there is no degree
preserving vector space isomorphism between the two spaces.

\b \no{\it 1.3.5.  Example.} If $q_{12}q_{21} = 1$, then $A
:= \xi_1\xi_2 - q_{21}\xi_2\xi_1$ is a constant of $\cb_q$,
and so is $A' := A\xi_3 - q_{31}q_{32}\xi_3A$.  Generically,
the space of constants of degree 1 in each generator
$\xi_1,\xi_2,\xi_3$ is ${\bf C}A'$.  But $\hat A :=
\sl_1\sl_2 - q_{21}\sl_2\sl_1$ belongs to Ker$D_q$ and
therefore so do $\hat A\sl_3$ and $\sl_3\hat A$.  So
Ker$D_q$ is bigger than the space of constants.  We shall
show that Ker$D_q$ is isomorphic to the ideal in $\cb_q$
generated by the constants.

\b \no{\bf 1.4.  More about constants.}
\b

The significance of the constants lies in their relation to
problems of integrability.

\b \no{\it 1.4.1.  Proposition.} Fix $q$ and a positive
integer $n$.  The following statements are equivalent:

\no (a) There are no constants in $\cb_q(n)$.

\no (b) The equations
$$
\partial_ix = y_i,\quad i = 1,...,N,
$$
have a solution $x \in \cb_q(n)$ for arbitrary $y_1,...,y_N
\in \cb_q(n-1)$.  In this case the solution is unique.

\b \no{\it Proof.} A straightforward extension of the proof
of Lemma 1.2.2.

\b \no{\it 1.4.2.  Proposition.} Let $C$ be a homogeneous
constant, and $x \in \cb_q$ any monomial of total degree
$k$.  There exists a constant of the form
$$
C' = \sum_{m = 0}^k\sum_{{\underline i}}
\xi_{i_1}...\xi_{i_m}a(\underline i) =
Cx + \sum_{m =
1}^k\sum_{{\underline i}}
\xi_{i_1}...\xi_{i_m}Ca(\underline i),
$$
where $a(\underline i)$ is a monomial of degree $k-m$.

\b \no{\it Proof.} There is a simple construction, using
induction on the degree $k$ of $x$, that leads to a
(generally unique) constant of this form.

\b \no{\it 1.4.3.  Definition.} Let ${\cal
I}_q$ denote the two-sided ideal of $\cb_q$ generated by the
constants in $\cb_q$.  Let $ \cb'_q$ be the quotient algebra
$\cb_q/{\cal I}_q$, and $\pi$ the projection of $\cb_q$ on
$\cb'_q$.  Since ${\cal I}_q$ is invariant under
differentiation, there is a natural action of $\partial_i$
on $\cb'_q$, namely $\partial_i \pi x = \pi \partial_i x$.
Thus $\cb'_q$ inherits the differential structure of
$\cb_q$.

\b \no{\it 1.4.4.} Let $ \widehat J $ be the
isomorphism $\cb \rightarrow \widehat{\cb}$ such that $\xi_i
\mapsto \sl_i$, and $J_q = D_q\circ\widehat J: \cb_q \rightarrow
\cb^*_q$ the unique homomorphism that maps $\xi_i$ to
$\partial_i$, $i = 1,...,N$.  \b \no{\bf 1.4.5.  Theorem.}
The mapping $\widehat{J}: \cb \rightarrow \widehat{\cb}$
induces an isomorphism $ {\cal I}_q \rightarrow {\rm Ker}
D_q$, and $J_q$ induces an isomorphism $\cb'_q \rightarrow
\cb^*_q$.  Hence $\cb^*_q$ is the topological dual of
$\cb'_q$.  \b The proof is in 1.6.  \b

\no
{\bf 1.5. The symmetric form S.}
\b

This will prepare the way for a proof of theorem 1.4.5.

\b \no{\it 1.5.1.  Definition.} Denote by $S_q$ the 2-form
on $\cb_q$ defined by \footnote*{This form was studied by
Kashiwara [K], and by Varchenko; both in the context of
Kac-Moody algebras.  It appears in the present, wider
context in a study of the standard universal R-matrix [F].}
$$
S_q(x,y) = \bigl((J_qx)y \bigr)_0,~~~x,y  \in \cb_q.
$$
Here $(\, )_0$ is the projection $\cb_q \rightarrow\cb_q(0)$.

\b \no{\it 1.5.2.  Proposition.} The form $S_q$ is
symmetric.

\b \no{\it Proof.} It is enough to consider the case when
$x$ and $y$ are monomials of the same degree.  Pairing the
operator $\partial_\alpha$ with each $\xi_r$ in turn we have
$$
\eqalign{
(...\partial_\gamma\partial_\beta\partial_\alpha
\xi_a\xi_b\xi_c...)_0 =
&\delta_{\alpha a}(...\partial_\gamma\partial_\beta
\xi_b\xi_c...)_0 +
\delta_{\alpha b}q_{\alpha a}(...\partial_\gamma
\partial_\beta\xi_a\xi_c...)_0\cr
& + \delta_{\alpha c}q_{\alpha a}q_{\alpha
b}(...\partial_\gamma\partial_\beta\xi_a\xi_b...)_0
+ ...~.
}
$$
A similar pairing of $\xi_\alpha$ with each $\partial_r$ gives
$$
\eqalign{
(\partial_a\partial_b\partial_c......\xi_\gamma\xi_\beta
\xi_\alpha )_0 =
&\delta_{\alpha a}( \partial_b\partial_c......\xi_\gamma
\xi_\beta)_0 +
\delta_{\alpha b}q_{ba}(
\partial_a\partial_c ......\xi_\gamma\xi_\beta )_0\cr
& + \delta_{\alpha c}q_{ca}q_{cb}(
\partial_a\partial_b......\xi_\gamma\xi_\beta)_0
+ ...~.
}
$$
The result follows by induction on the degree.
\bb

\no {\bf 1.6.  Proof of the theorem 1.4.5.}

\b \no{\it 1.6.1.} Let $x \in {\cal I}_q$, a sum of
homogeneous polynomials each of which contains a constant
factor.  It was shown that $\widehat J$ maps constants into
Ker$D_q$; it defines a monomorphism from ${\cal I}_q$ into
Ker$D_q$.

\b \no{\it 1.6.2.} Conversely, let $x'\in \widehat {\cb}$,
then there is $x \in \cb$ such that $x' = \widehat {J}x$.
In particular, let $\widehat{J}x \in {\rm Ker}D_q$,
homogeneous of degree $k$.  Then for all $y \in \cb_q$ of
the same degree, $(J_qx)y = 0$, and by the symmetry of $S$,
$$
(J_qy)x = 0,\quad \forall~ \widehat{J}x \in {\rm Ker}\,
D_q,~{\rm deg}(y) =
{\rm deg}(x).
$$
It remains to be shown that this result implies that $x$
belongs to ${\cal I}_q$.

\b \no{\it 1.6.3.  Lemma.} (a) There are no constants in
$\cb_q'$.

\no (b) If $x \in \cb_q$ is
homogeneous of degree $k$, and if $x$ is annihilated by all differential
operators of
degree $k$, then $x \in {\cal I}_q$.

\no (c) Every $x \in \cb$ can be expressed as
a ``Taylor series",
$$
 x = c(x) +
 \sum_{n \geq 1}\sum_{{\underline i}}A^{{\underline i}}
 \partial_{i_1}...
\partial_{i_n}x,\quad
 A^{{\underline i}} = \sum_\sigma A^\sigma
 \xi_{i_{\sigma 1}}...\xi_{i_{\sigma n}},
$$
where $\partial_ic(x)=0, i=1,...,N$; the second sum is over
the permutations of $1,...,n$, and $A^\sigma \in {\bf C}$.
(The coefficients are universal, independent of $x$; they
are calculated in 1.6.5 and examples are given in 1.6-7.)

\b \no{\it Proof.} (a) means that, if $\partial_i \pi x =
0$, then $\pi x = 0$; that is, if $\pi$ annihilates the
derivatives of $x$, then $\pi$ annihilates $x$.  Suppose
that $x$ is homogeneous of degree $k$, and that all
derivatives of order $k$ are zero,
$\partial_{i_1}...\partial_{i_k} x = 0$.  Then
$\pi\partial_{i_1}...\partial_{i_k} x =
\partial_{i_1}...\partial_{i_k}\,\pi x = 0$.  Assume (a),
then it follows that $\pi x = 0$.  Hence (a) implies (b).
On the other hand, (c) implies (a) (as is seen by applying
$\pi$ to both sides of the formula) so it is enough to prove
(c).

To prove (c), it is enough to prove that the coefficients
$A^{{\underline i}}$ can be chosen so that both sides of the
equation have the same first derivatives, which means that
$$
\sum_{{\underline i}}\bigl(
\partial_kA^{ i_1...i_n} + \delta_{k,i_1}q_{i_1i_2}...q_{i_1i_n}
A^{i_2...i_n} \bigr)\partial_{i_1}...\partial_{i_n}  =
0,\quad n = 1,2,...~
.\eqno(1.1)
$$
Here $A^{i_1...i_n}$ is defined to be equal to $-1$ for $n =
0$.  This looks like a sequence of equations of the form
$\partial_ix = y_i$, to the solutions of which the presence
of constants represents an obstruction.  But the
obstructions are, in fact, circumvented.  Consider the
operator that sends $A\in \cb(n)$ to the $\cb(1)$ one-form
$(\partial_kA)_{k = 1,...,N}$, valued in $\cb(n-1)$.  On
suitable bases, denote by $M$ the associated square matrix.
The obstructions to solving
$$
\partial_kA^{ i_1...i_n} + \delta_{k,i_1}q_{i_1i_2}...q_{i_1i_n}
A^{i_2...i_n} = 0,\eqno(1.2)
$$
are the constants in $\cb(n)$, the null vectors of $M$.  To
the null space of $M$ there corresponds the null space of
the transposed matrix, of the same dimension.  Indeed, the
map $J_q$ defined in 1.4.4 induces a natural bijection from
one to the other.  Let
$$
C = \sum_{\underline j} C^{j_1...j_n}\xi_{j_1}...\xi_{j_n} =
\sum_{\underline j}C^{\underline j}\,
\xi_{j_1}...\xi_{j_n}
$$
be a constant. Then
$$
J_qC = \sum_{\underline j} C^{\underline j}\,\partial_{j_1}...
\partial_{j_n}
$$
vanishes identically and this represents a null vector of
the transposed of $M$:
$$
\sum_{\underline j}\bigl(C^{\underline j}\,\partial_{j_1}
...\partial_{j_{n-1}}\bigr)\partial_{j_n}A = 0.\eqno(1.3)
$$
Now assume that Eq.(1.1) can be solved for $n = 1,...,m-1$,
then (1.2) is valid for $n = 1,...,m-1$; not identically, but
as a substitution under
$\sum\partial_{i_1}...\partial_{i_n}$.  Eq.(1.3) says that the
one-form $(\partial_kA)$ is closed; the obstruction to
solving Eq.(1.2) consists of the fact that the second term is
not closed.  However, using (1.2) for $n = 1,...,m-1$, in the
sense just explained, we find that
$$
\sum_{\underline j}\bigl(C^{\underline j}\,\partial_{j_1}...
\partial_{j_{n-1}}\bigr) \delta_{j_n,i_1}q_{i_1i_2}
...q_{i_1i_n}A^{i_2...i_n} \propto C^{{\underline i}}
\approx 0.
$$
Therefore, the obstructions to solving Eq.(1.2) do not affect
Eq.(1.1), and we conclude that (1.1) is always solvable, giving
a unique solution for
$\sum_{\underline i}A^{i_1...i_n}\partial_{i_1}...\partial_{i_n}$.  The lemma
is established by induction in $n$, and with that, Theorem
1.4.5 is proved.

\b \no{\it 1.6.4.  Corollary.} The radical of the form $S_q$
is the ideal ${\cal I}_q$.  By the projection $\pi: \cb_q
\rightarrow \cb'_q$, we get a nondegenerate two-form on
$\cb'_q$ that will also be denoted $S_q$; it can be
interpreted as an invertible map $S_q: \cb'_q \rightarrow
\cb^*_q$.

\b \no{\it 1.6.5.} We can actually determine the
coefficients $A^{{\underline i}}$ explicitly.  To this end
let $q$ be in general position and set
$$
A^{i_1...i_n} = (-)^{n+1}\bigl(\prod_{k<l} q_{i_ki_l}\bigr)
T^{i_1...i_n},
$$
 Then the zero-rank tensor $T = 1$, and the recursion relation (2) reduces to
$$
\partial_kT^{i_1...i_n} = \delta_{k,i_1}T^{i_2...i_n}.
$$
Iteration gives
$$
\partial_{k_n}...\partial_{k_1}T^{{\underline i}}=
\delta_{k_1,i_1}...\delta_{k_n,i_n}.
$$
Setting
$$
T^{{\underline i}} =
\sum_{{\underline j}}T^{{\underline i}}_{{\underline j}}\,
\xi_{j_1}...\xi_{j_n}
$$
we obtain
$$
\sum_{{\underline j}} T^{{\underline i}}_{{\underline j}}\,
S_{{\underline k}}^{{\underline j}} = \delta_{k_1,i_1}
...\delta_{k_1,i_n},
$$
where the coefficients of the form $S$ are defined by
$$
S_{i_1...i_n}^{j_1...j_n}=
\partial_{i_n}...\partial_{i_1}(\xi_{j_1}...\xi_{j_n}).
$$
Hence $T$ can be interpreted as the contragredient two-form,
inverse to $S$.  This interpretation survives at exceptional
points in the space of parameters where $S$ is a
non-degenerate two-form on $\cb_q'$.  The coefficients
$T^{{\underline i}}_{{\underline j}}$ are not unique, but
$\sum_{{\underline j}}T^{{\underline i}}_{{\underline j}}\,
\partial_{j_1}...\partial_{j_n}$ is unique and so is $\pi
\,\sum_{{\underline i}}T^{{\underline i}}_{{\underline j}}\,
\xi_{j_1}...\xi_{j_n}$.

\bb
\noindent {\it 1.6.6.  Example.} Suppose $N=1$.  If $q:=q_{11}$
is not a root of unity, then, for all $x\in {\cal B}_q$,
$$
x = x_0+\xi_1\partial_1 x + \sum_{n\geq 2}(-1)^{n-1}
{q^{n\choose 2}\over
[n]_q!}\xi^n_1\partial^n_1x,\quad [n]_q!:=
\prod_{k=1}^n[k]_q,\quad [k]_q:=\sum_{j=0}^{k-1}q^j,
$$
where $x_0$ is the projection of $x$ on ${\cal B}(0)$.  When
$q^m=1, q\neq 1$, the expansion truncates at $n=m-1$, and
$x_0$ is replaced by $c(x)$, with $\partial_1 c(x)=0$.

\bigskip
\noindent{\it 1.6.7.  Example.} Suppose $N$ arbitrary.  If
$q_{ij}$ are in general position, then, for all $x\in {\cal
B}$,
$$
\eqalign{
x &=  x_0 + \sum_{i}\xi_i\partial_ix-\sum_i {q_{ii}\over
[2]_{q_{ii}}!}\xi_i^2\partial_i^2x
-\sum_{i\neq j}{q_{ij}\over
1-\sigma_{ij}}(\xi_i\xi_j-q_{ji}\xi_j\xi_i)
\partial_i\partial_jx  + \ldots,
}
$$
where $x_0$ is the projection of $x$ on ${\cal B}(0)$.  When
there is the  constraint $\sigma_{12}=1$, the above
expansion holds provided that $x_0$ is replaced by some
$c(x)$, with $\partial_i c(x)=0, i=1,\ldots,N$.
Terms that blow up due to the factor $(1-\sigma_{12})$ in
the  denominators are  replaced by
expressions that can be obtained either  solving (1)
 or from the generic expansion by a limiting
procedure; thus
$$
\eqalign{
&-\sum_{i\neq j=1,2}{q_{ij}\over
1-\sigma_{ij}}(\xi_i\xi_j-q_{ji}\xi_j\xi_i)
\partial_i\partial_jx \cr
&\mapsto
(\alpha q_{12}\xi_1\xi_2+\beta \xi_2\xi_1)\partial_1\partial_2x
+(\gamma \xi_1\xi_2+\delta q_{21}\xi_2\xi_1)\partial_1\partial_2x
}
$$
where $\alpha,\beta,\gamma,\delta\in{\bf C}$ are arbitrary
complex numbers satisfying $\alpha+\beta+\gamma+\delta+1=0$.

Ambiguities in these parameters are absorbed into the
constant term.

\bb
\no{\steptwo 2. Yang-Baxter algebras.}
\b

Our interest now focuses (in this section only) on the
exceptional values of the parameters and on the quotient
algebras $\cb'_q$.  Here we shall connect all these algebras
to{\it Yang-Baxter} algebras and some of them to quantized
Kac-Moody algebras.

To every $\cb'_q$ corresponds a Yang-Baxter algebra.  This
family of algebras includes the quantized Kac-Moody algebras
and complete them in a natural way.  They all admit a
coboundary Hopf structure with a universal R-matrix of a
standard form.  (The Hopf structure is reviewed in [F].)

\b \no{\bf 2.1.  Generators.} \b The generators
($\xi_i,\partial_i)$ of $\cb_q$ are related to the
Chevalley-Serre generators $(e_i,f_i)$ of quantized
Kac-Moody algebras.  We first introduce the ``Cartan
subalgebra".  \b \no{\it 2.1.1.} Let $\hat K^i$ and $\hat
K_i$ be the unique automorphisms of $\cb_q$ such that
$$
\hat K^i(\xi_j) = q_{ij}\xi_j,\quad
\hat K_i(\xi_j) = {1\over q_{ji}}\xi_j.
$$
We expand the algebra $\cb_q$ by including new generators
$K^i,K_i, i = 1,...,N$ that implement these automorphisms,
thus
$$
  K^i \xi_j  = q_{ij}\xi_jK^i,\quad
  K_i \xi_j  = {1\over q_{ji}}\xi_jK_i,\quad i,j = 1,...,N.
$$

\b \no{\it 2.1.2.} Let $(\overleftarrow\partial_i)_{ i =
1,...,N}$ be differential operators that act on $x \in
\cb_q$ from the right, such that
$$
x\xi_i\overleftarrow\partial_j = x\delta_{ij} +
 xq_{ij}\overleftarrow\partial_j\xi_i.
$$

\b \no{\it 2.1.3.  Proposition.} [F] The ideal in $\cb_q$
that is generated by the constants with respect to the
operators $\partial_i$ coincides with the ideal generated by
the constants with respect to the operators
$\overleftarrow\partial_i$.

\b \no{\it 2.1.4.} For any $q = \{q_{ij}\}$, set
$$
e_i = \xi_i,\quad f_i =
\overleftarrow\partial_iK_i - K^i\partial_i,
$$
then the following relations hold for $i,j = 1,...,N$:
$$
\eqalign{
&[K_i,K_j] = [K_i,K^j] = [K^i,K^j] = 0,\cr
& K^i e_j  =  q_{ij}e_jK^i,~~~   K_ie_j  = (q_{ji})^{-1}e_jK_i,\cr
& K^i f_j  =  (q_{ij})^{-1}f_jK^i,\quad   K_if_j  = q_{ji}f_jK_i,\cr
&[e_i,f_j] = \delta_{ij}(K^i - K_i),\quad i,j = 1,...,N.
}
$$
These are the relations of (the multiparameter version of)
Drinfel'd's quantization of Kac-Moody algebras, except (i)
for the omission of Serre relations and (ii) certain
conditions on $q$ that we discuss next.  Of course, (i) and
(ii) are very closely connected.

\b \no {\bf 2.2.  Serre relations.} \b The ideals ${\cal
I}_q$ of $\cb_q$ that appear for exceptional values of the
parameters have not yet been classified, but it is not
difficult to find examples.  We specialize, temporarily, to
quantized Kac-Moody algebras.

Suppose that for each ordered pair $(i,j)_{i,j = 1,...,N}$ there is a
positive integer
$k_{ij}$ such that
$$
\sigma_{ij}q_{ii}^{k_{ij}-1} = 1,  ~~\sigma_{ij} := q_{ij}q_{ji}.
$$
In this case the following elements of $\cb_q$ are constants,
$$\eqalign{
&\sum_{m=0}^{k_{ij}}Q_{km}^{ij}(\xi_i)^m\xi_j(\xi_i)^{k_{ij}-m},\cr &
Q_{km}^{ij} :=  (-q_{ij})^m(q_{ii})^{m(m-1)/2}
{k_{ij}\choose m}_{q_{ii}},\quad
{k\choose m}_{q } := {[k]_q!\over [m]_q![k-m]_q!},
}
$$
and the image of each one by $J_q$ is identically zero.  The
correspondence in 2.1.4, and passage to the quotient, now
yields the quantized Kac-Moody algebra with Cartan matrix
$$
A_{ij} = 1 - k_{ij}.
$$

\ve
\b \no{\bf 2.3.  An alternative presentation.}

\b \no{\it 2.3.1.  Definition.} Let ${\cal M,N}$ be two
finite sets, $\varphi,\psi$ two maps,
$$
\eqalign{
& \varphi :~~{\cal M}\otimes {\cal M}
\rightarrow {\bf C}, \cr & \psi
:~~{\cal M}\otimes {\cal N} \rightarrow {\bf C}, \cr} \quad
\eqalign{a,b &\rightarrow \varphi^{ab}~, \cr a,i & \rightarrow
H_a(i)~. \cr}
$$
Let ${\cal{A}}$ or ${\cal{A}}(\varphi,\psi)$ be the
universal, associative, unital algebra over {\bf C} with
generators $\{H_a\}_{a\in {\cal M}}$, $\{e_i, f_i\}_{i \in
{\cal N}}$, and relations
$$
\eqalign{&  [H_a,H_b]=0~,\cr &  [H_a,e_i] = H_a(i)
e_i,\quad [H_a,f_i]
=-H_a(i)f_i,\cr
  &[e_i,f_j]=\delta_{ij}
\bigl({\rm e}^{\varphi(i,\cdot)}-{\rm e}^{-\varphi(\cdot,i)}
\bigr)~,
  \cr}
$$
with $\varphi(i,\cdot)=\sum_{a,b}\varphi^{ab}H_a(i)H_b,~
\varphi(\cdot,i)= \sum_{a,b}\varphi^{ab}H_aH_b(i),~
{\varphi(i, \cdot) + \varphi(\cdot,i)} \neq 0,~ i \in {\cal
N}$.  (The last condition on the parameters is included in
order to avoid having to make some rather trivial
exceptions.)  For $H\in {\cal A},~ {\rm e}^H $ is the formal
series ${\rm e}^H := \sum_n{1\over n!}H^n $.  We take ${\cal
N} = \{1,...,N\}$.

These relations imply those in 2.1.4 if we set
$$
K^i = {\rm e}^{\varphi(i,.)},\quad K_i = {\rm e}^{-\varphi(.,i)},
$$
and choose $\varphi$ so that
$$
{\rm e}^{\varphi(i,j)} = q_{ij},\quad \varphi(i,j) :=
\sum_{a,b}\varphi^{ab}H_a(i)H_b(j).
$$
This alternative presentation is more cumbersome, but it
seems to be necessary for the introduction of a Universal
R-matrix.

\b \no{\it 2.3.2.} The free subalgebra generated by
$\{e_i\}_{i = 1,...,N}$ (resp.  $\{f_i\}_{i= 1,...,N}$) will
be denoted ${\cal{A}}^+$ (resp.  ${\cal{A}}^-$).  The {\it
Cartan subalgebra} generated by $\{H_a\}_{a\in {\cal M}}$ is
denoted ${\cal{A}}^0$.  The passage to the quotient $\cb'_q$
gives a quotient ${\cal A}'$ with subalgebras
${\cal{A}}'^\pm, {\cal{A}}'^0$.

\b \no {\it 2.3.3.  Yang-Baxter element. (Universal {\rm
R}-matrix.)}

The Universal R-matrix for ${\cal A}$ exists for parameters in
general position and is given by [F]
$$
R = {\rm e}^\varphi\sum_{n=0}^\infty t_n,\quad
t_n = \sum_{{\underline i},{\underline j}}
T_{{\underline i}}^{{\underline j}}f_{i_1}...f_{i_n}
\otimes e_{j_1}...e_{j_n}.
$$
The ``Cartan potential",
$${\rm e}^\varphi := {\rm e}^{\sum_{a,b}\varphi^{ab}
H_a\otimes H_b},
$$
satisfies the following relations
$$
{\rm e}^\varphi(e_i \otimes 1) =
(e_i \otimes K^i){\rm e}^\varphi,\quad
 (1 \otimes e_i){\rm e}^\varphi=
 {\rm e}^\varphi(K_i \otimes e_i).
$$
For exceptional values of the parameters the expression for
$R$ makes sense in ${\cal A}' \otimes {\cal A}'$, and gives
the standard Universal R-matrix for the quotient algebra
${\cal A}'$.

\b \no {\it 2.3.4.  Remark.} For every set of parameters
$q_{ij}, q_{ii} \neq 1, 1 = 1,...,N$, there are $A_{ij} \in
{\bf C}$ such that $\sigma_{ij} = q_{ii}^{A_{ij}}, i\neq j$,
$A_{ii} = 2$, and such that the matrix $(A_{ij})$ is
symmetrizable.  Any such matrix may be called the Cartan
matrix of $\cb_q$.

\bb
\no{\bf 2.4.  Highest root vectors and quantized loop
algebras.}
\b
\no{\it 2.4.1.} We ask whether the quantized
enveloping algebra of a finite dimensional Lie algebra, of
standard (multidimensional) type, admits a highest root
generator; that is, an element $E\in U_q($\g$)$, with the
same weight relatively to the Cartan subalgebra as the
highest root in \g, and a set of complex coefficients
$a_1,...,a_l$ ($l$ is the rank of \g) such that
$$
w_i:=[e_i,E]_{a_i} = 0,~~ i = 1,...,l.\eqno(2.1)
$$
An affirmative answer is given for \g ~in the series $A_l$
and $C_l$; and a negative one for $B_3, G_2$.  Proofs are
included for the case of $C_l$ only.  When $E$ exists, let
$F$ be the corresponding negative root generator.  Let
$U_q^z($\g$) = U_q($\g$)\otimes{\bf C}[z,1/z]$; this
algebra, generated by the set of generators of $U_q($\g$)$
augmented by $e_0 := zF$ and $f_0 := z^{-1}E$, satisfies
all the relations of a quantized Kac-Moody algebra of
affine type except, possibly, for the relation $[e_0,f_0] =
{\rm e}^{H(0)} - {\rm e}^{-H(0)}$.

\b \no{\it 2.4.2.  The case of $A_l$.}

The constraints on the parameters are
$$
\eqalign{&\sigma_{ij} = 1,~~ |i-j| > 1,
~~\sigma_{ij} := q_{ij}q_{ji},~
i\not= j, \cr
&\sigma_{ij}q_{jj} = 1,~~ |i-j|=1,
~~i,j = 1,\ldots,\ell,\cr} \eqno(2.2)
$$
\no and the ideal ${\cal I}_q$ is generated by
$$
\eqalign{&\xi_i\xi_j-q_{ji}\xi_j\xi_i, \quad |i-j| > 1,
\quad i,j = 1,\ldots,\ell, \cr
&C_{iij} := {1\over q_{ji}} \xi_i^2\xi_j-
[2]_{q_{ii}}\xi_i\xi_j\xi_i +
{1\over q_{ij}}\xi_j\xi_i^2, ~~ |i-j|=1. \cr} \eqno(2.3)
$$
Following Jimbo [J], we introduce a sequence of ``root vectors"
$$
\eqalign{
&X^1 := \xi_1, ~~X^n = [X^{n-1},\xi_n]_{a_n},~~
n=2,\ldots,\ell, \cr
&a_n = q_{n1}\ldots q_{n,n-2}/q_{n-1,n}. \cr} \eqno(2.4)
$$

\b\no {\bf Proposition 2.4.3.}  Let $X^0=1$, then
$$
\eqalign{&\partial_iX^n \propto \delta_{in} X^{n-1},
\quad ~~ i,n=1,\ldots,\ell, \cr
&[X^n,\xi_i]_{q_{i1}\ldots q_{in}}=0, \quad i=n+2,\ldots,\ell. \cr} \eqno(2.5)
$$
\no The constants of proportionality are $\not= 0$.

\b\no {\bf Proposition 2.4.4.}  With
$$
\eqalign{&a(i,n) = q_{i1}\ldots q_{in}, \quad i \not= n,n+1, \cr
&a(n,n) = q_{n1}\ldots q_{n,n-1}, \quad n\not= 1, ~a(1,1)=1,\cr
&a(n,n-1) = a_n = q_{n1}\ldots q_{n,n-2}/q_{n-1,n}~,\cr} \eqno(2.6)
$$
one has for $n,i=1,\ldots,\ell$,
$$
[X^n,\xi_i]_{a(i,n)} = \delta_{i,n+1}X^{n+1}. \eqno(2.7)
$$

\b\no {\bf Corollary 2.4.5.} The element $X^\ell$ is a
highest root vector for $A_\ell$.

\b \no {\it 2.4.6.  The case of $B_2, B_3$.}

Here, we had intended to present our study of $B_l$, but
when this case turned out to be more difficult we scaled
back our ambition and attacked $B_2$ and $B_3$ instead.  In
the first case there is a unique (up to normalization)
highest weight generator.  The unexpected negative result
for $B_3$ is stated as Theorem 2.4.8, below.

\b We begin with $B_2$.  The data encoded in the Cartan
matrix is $ q_{11} = q_{22}^2 = 1/\sigma_{12}$.  The long
simple root corresponds to $\xi_1$ and $E \in \cb_{122}$;
that is, $E$ is of first order in $\xi_1$ and of second
order in $\xi_2$.  The ideal ${\cal I}_q$ is generated by
$C_{112}$, Eq. (2.3), and
$$\eqalign{
&C_{2221} := {1\over q_{12}}\xi_2^3\xi_1
- [3]_{q_{22}}\xi_2^2\xi_1\xi_2 +
[3]_{q_{22}}q_{12}q_{22}\xi_2\xi_1\xi_2^2 -
q_{12}^2q_{22}^3\xi_1\xi_2^3.\cr}\eqno(2.8)
$$
As the space of constants in $\cb_{1222}$ is 1-dimensional,
we may set $w_2:=[E,\xi_2]_{a_2}\propto C_{2221}$; that is,
we look for an element
$E$ such that
$$
C_{2221} \propto E\xi_2 - a_2\xi_2E.
$$
There are exactly 3 (linearly independent) solutions of this
equation.  They are
$$
\eqalign{ &E = \xi_1\xi_2^2 -
q_{21}q_{22}(1+q_{22})\xi_2\xi_1\xi_2 +
q_{21}^2q_{22}^3\xi_2^2\xi_1, ~a_2 = q_{21},\cr &E =
\xi_1\xi_2^2 - q_{21}(1+q_{11})\xi_2\xi_1\xi_2 +
q_{21}^2q_{22}^2\xi_2^2\xi_1, ~~~~~a_2 = q_{21}q_{22},\cr &E
= \xi_1\xi_2^2 - q_{21}(1+q_{22})\xi_2\xi_1\xi_2 +
q_{21}^2q_{22}\xi_2^2\xi_1,~~~~ ~a_2 = q_{21}q_{22}^2.\cr}
$$

The task is finished if there is an element in ${\cal I}_q$
of the form $w_1:= E\xi_1 - a_1\xi_1E$.  It is easy to see
that $\partial_1w_1 = 0$ requires that $\partial_1E = 0$.
The first two solutions satisfy this requirement.  There
remains to satisfy $\partial_2w_1 = 0$, and this makes the
solution unique:
$$
E = C_{221} = [[\xi_1,\xi_2]_{q_{21}},\xi_2]_{q_{21}q_{11}},
~~ a_1 =
q_{12}^2q_{22}^2,~~ a_2 = q_{21}q_{22}.
$$
Note that $C_{221}$ is not a constant in this case.

\b
On to $B_3$.
\b

\no{\it Remark 2.4.7.} A convenient basis for $B_3 \sim
so(7)$ may be found in [J].  If $E_{ij}$ is the $7\times 7$
matrix that has 1 in position $(i,j)$ and 0's elsewhere,
then a system of simple root generators is
$$
e_1 = E_{12}-E_{67},~~e_2 = E_{23}-E_{56},~~e_3 = E_{34}-E_{45}.
$$
The highest root is $r_1 + 2r_2 + 2r_3$, and $r_3$ is the
short root.  The highest root generator can be expressed in
several different ways, but if it is written as $[A,e_i], A
\in B_3$, then always $i = 2$,
$$
E = [A,e_2].\eqno(2.9)
$$
\b In the quantized case, we would need an element $E \in
\cb'_{12233}$; that is, a sum of permutations of
$\xi_1\xi_2\xi_2\xi_3\xi_3$, and constants $a_1,a_2,a_3$,
such that
$$
w_i := [E,\xi_i]_{a_i} = 0,~~ i  = 1,2,3.\eqno(2.10)
$$
(We do not assume that $E$ has the form $[A,\xi_2]_a$
suggested by (2.9).)

\b \no{\bf Theorem 2.4.8.} There is no element $E\neq 0$ in
$\cb'_{12233}$ that satisfies (2.10), and that tends
to the highest root generator of $B_3$ in the Lie limit.  \b
\no {\it 2.4.9.  The case of $C_\ell$.}

The highest root of $C_\ell$, in terms of the simple roots,
is $2r_1+\ldots + 2r_{\ell-1}+r_\ell$; where $r_\ell$ is the
long root.  To the sequence
$$
1, 2, ... , (l-1) , l, (l-1) , ... , 2, 1\eqno(2.11)
$$
we associate a sequence $X^1,...,X^l,X^{l+1},...,X^{2l-1}$
of root vectors defined as follows,
$$
\eqalignno{
&X^1 = \xi_1, \quad X^n = [X^{n-1},\xi_n]_{a_n},
\quad n=2,\ldots,\ell, & (2.12) \cr
&X^{\ell+n} = [X^{\ell+n-1},\xi_{\ell-n}]_{b_n}, \quad
~\,n=1,\ldots,\ell-1, & (2.13) \cr}
$$
with coefficients $a_n,b_n$ to be chosen.

The constraints are
$$
\eqalign{&\sigma_{ij}=1,~~|i-j|>1, \cr
&\sigma_{ij}q_{ii}=1,~~|i-j|=1,~~i,j=1,\ldots,\ell-1, \cr
&\sigma_{\ell,\ell-1}q_{\ell\ell}=1,~~
q_{11}=\ldots=q_{\ell-1,\ell-1} =: q,~~
q_{\ell\ell}=q^2. \cr} \eqno(2.14)
$$
The ideal is generated by (2.2), except that
$C_{\ell-1,\ell-1,\ell}$ is replaced by
$C_{\ell-1,\ell-1,\ell-1,\ell}$ defined in (2.8).

We take over the definition (2.4) of $X^1,\ldots,X^\ell$.,
and the relations (2.5) still hold.  Of the relations (2.7),
all but one remains valid;
$[X^\ell,\xi_{\ell-1}]_{a(\ell-1,\ell)}$ is no longer zero
because the constraint
$\sigma_{\ell,\ell-1}q_{\ell-1,\ell-1}=1$ no longer applies,
since $r_{l-1}$ is short.  Thus
$$
[X^n,\xi_i]_{a(i,n)} = \delta_{i,n+1}X^{n+1},
~~\matrix{ n,i=1,...,l,\cr{\rm
except ~for~} n = l,~~i = l-1.\cr}
\eqno(2.15)
$$
This invites us to construct the sequence (2.13), the
properties of which we shall now explore.

\b\no {\bf Proposition 2.4.10.} For $n=1,\ldots,\ell-1$,
define
$$
\eqalign{
X^{\ell+n} &= [X^{\ell+n-1},\xi_{\ell-n}]_{b_n}, \cr
b_n &= q_{\ell-n,1}\ldots q_{\ell-n,\ell}q_{\ell-n,\ell-1} \ldots
q_{\ell-n,\ell-n} =: q_{\ell-n,1} \ldots,\ldots
q_{\ell-n,\ell-n}. \cr} \eqno(2.16)
$$
Then
$$
\partial_iX^{\ell+n} \propto \delta_{i,\ell-n}
X^{\ell+n-1}, n = 1,...,l-1,~~i = 1,...,l. \eqno(2.17)
$$

\b\no {\bf Proof.} For $n=1$ we have
$$
\eqalign{
\partial_iX^{\ell+1} &= \partial_i\,[X^\ell,
\xi_{\ell-1}]_{b_1}
= [\partial_iX^\ell,\xi_{\ell-1}]_{b_1q_{i,\ell-1}}+
\delta_i^{\ell-1}(q_{i1}\ldots q_{i\ell}-b_1)X^\ell \cr}
$$
This is zero for $i=1,\ldots,\ell-2$.  With the help of
(2.5) one gets
$$
\eqalign{
\partial_lX^{\ell+1} &\propto [X^{\ell-1},\xi_{\ell-1}]_
{b_1q_{\ell,\ell-1}} \cr
}
$$
The constraint $q^2_{\ell-1,\ell-1}q_{\ell-1,\ell}
q_{\ell,\ell-1}=1$ makes
$$
b_1q_{\ell,\ell-1} = q_{\ell-1,1}\ldots q_{\ell-1,\ell-2} =
a(\ell-1,\ell-1),
$$
and
$$
\eqalign{
\partial_lX^{\ell+1} &\propto
[X^{\ell-1},\xi_{\ell-1}]_{a(\ell-1,\ell-1)} = 0, \cr}
$$
Finally, $\partial_{\ell-1}X^{\ell+1}\propto X^\ell$, as
required.  This establishes a base for induction in $n$.
Suppose that the statement of the proposition is true for
$n=1,\ldots,m$; we have for $n = 1,...,l-1,\break~ i =
1,...,l$,
$$
\eqalign{\partial_iX^{\ell+m+1} &=
\partial_i[X^{\ell+m},\xi_{\ell-m-1}]_{b_{m+1}} \cr
&= [\partial_iX^{\ell+m},\xi_{\ell-m-1}]_{b_{m+1}}q_{i,\ell-m+1}
 +\delta_i^{\ell-m-1}
(q_{i1}...,... q_{i,\ell-m}-b_{m+1})X^{\ell+m}. \cr}
$$
By hypothesis this is zero for $i\not= \ell-m-1,~\ell-m$,
while
$$
\eqalign{\partial_{\ell-m}X^{\ell+m+1} &=
[\partial_{\ell-m}X^{\ell+m},\xi_{\ell-m-1}]_{b_{m+1}
q_{\ell-m,\ell-m-1}} \cr
&\propto [X^{\ell+m-1},\xi_{\ell-m-1}]_{b_{m+1}q_{\ell-m,
\ell-m-1}} \cr} \eqno(2.18)
$$
For $m=0$ this was already seen to vanish.  For $m=1$ it is
zero because of (2.15) and
$$
b_2q_{\ell-1,\ell-2}=q_{\ell-2,1}... q_{\ell-2,\ell} =
a(\ell-2,\ell).
$$
For $m>1$ we again invoke the principle according to which
(2.18) vanishes if and only all its derivatives vanish.  In
fact,
$$
\partial_i[X^{\ell+m-1},\xi_{\ell-m-1}]_{b_{m+1}
q_{\ell-m,\ell-m-1}}\eqno(2.19)
$$
vanishes for $i\not= \ell-m\pm 1$ by the induction
hypothesis, so we have to prove that it vanishes for
$i=\ell-m\pm 1$ as well.  First, taking $i=\ell-m-1$,
$$
\eqalign{
&\partial_{\ell-m-1}[X^{\ell+m-1},\xi_{\ell-m-1}]_
{b_{m+1}q_{\ell-m,\ell-m-1}} \cr
&= (q_{\ell-m-1,1}\ldots \ldots q_{\ell-m-1,\ell-m+1}-
b_{m+1}q_{\ell-m,\ell-m-1})X^{\ell+m-1} = 0, \cr}
$$
by virtue of the constraint
$\sigma_{\ell-m,\ell-m-1}q_{\ell-m-1,\ell-m-1}=1,~~ m>0$,.
For $i=\ell-m+1$, (2.19) becomes
$$
\eqalign{
&\partial_{\ell-m+1}[X^{\ell+m-1},\xi_{\ell-m-1}]_
{b_{m+1}q_{\ell-m,\ell-m-1}} \cr
&\propto [X^{\ell+m-2},\xi_{\ell-m-1}]_
{b_{m+1}q_{\ell-m,\ell-m-1}q_{\ell-m+1,\ell-m-1}} \cr}
\eqno(2.20)
$$
For $m=2$ it is
$$
\eqalign{
&[X^\ell,\xi_{\ell-3}]_{q_{\ell-3,1}...,...
q_{\ell-3,\ell}q_{\ell-3,\ell-1}q_{\ell-3,\ell-2}
q_{\ell-3,\ell-3}/q_{\ell-3,\ell-2}q_{\ell-3,\ell-3}
q_{\ell-3,\ell-1}} \cr
&= [X^\ell,\xi_{\ell-3}]_{q_{\ell-3,1}\ldots q_{\ell-3,\ell}} =
[X^\ell,\xi_{\ell-3}]_{a(\ell-3,\ell)} = 0. \cr}
$$
Now (2.20) looks like a shifted form of (2.18) and this
suggests to repeat the stops that led from one to the other.
Thus, to show that (2.20) is zero, it is enough to verify
that the expression is annihilated by $\partial_{\ell-m+2}$
and by $\partial_{\ell-m-1}$.  That $\partial_{\ell-m-1}$
gives zero is obvious since this operator quommutes with
$\partial_{\ell-m+1}$.  (Two operators $A,B$ quommute if
there is $\alpha$ in the field such that $[A,B]_\alpha := AB
- \alpha BA = 0$.)  So it is enough to consider
$$
\eqalign{
&\partial_{\ell-m+2}
[X^{\ell+m-2},\xi_{\ell-m-1}]_{b_{m+1}q_{\ell-m,\ell-m-1}
q_{\ell-m+1,\ell-m-1}} \cr
&\propto [X^{\ell+m-3},\xi_{\ell-m-1}]_{b_{m+1}
q_{\ell-m,\ell-m-1}q_{\ell-m+1,\ell-m-1}
q_{\ell-m+2,\ell-m-1}}. \cr}
$$
\no This vanishes for $m=3$ since the coefficient is then
$$
\eqalign{
&q_{\ell-4,1}\ldots q_{\ell-4,\ell}
q_{\ell-4,\ell-1} \ldots
q_{\ell-4,\ell-4}/q_{\ell-4,\ell-4}q_{\ell-4,\ell-3}
q_{\ell-4,\ell-2}q_{\ell-4,\
ell-1} \cr
&= a(\ell-4,\ell). \cr}
$$
The pattern is clear, after $m-1$ iterations we end up with
$[X^\ell,\xi_{\ell-m-1}]_{a(\ell,\ell-m-1)} = 0$.
The theorem is proved.

\b We collect all the relations obtained so far, and some
new ones.

\b \no {\it 2.4.11.  Properties of root vectors.} (a) From
(2.15) and (2.16), now valid for $i,n=1,\ldots,\ell$ except
for $n=\ell,i=\ell-1$:
$$
[X^n,\xi_i]_{a(i,n)} = \delta_i^{n+1}X^{n+1}, \eqno(2.21)
$$
where $a(i,n)$ is as in (2,6).

(b) From (2.13),
$$
[X^{\ell+n-1},\xi_{\ell-n}]_{b_n} = X^{\ell+n}, \quad
n=1,\ldots,\ell-1, \eqno(2.22)
$$
where $b_n$ is as  in (2.16).

(c) For $m=0,\ldots,\ell-1,~
n=2,3,\ldots,2m$,
$$
[X^{\ell+m-n},\xi_{\ell-m}]_{b(m,n)} = 0, \eqno(2.23)
$$
with $b(m,n) = q_{\ell-m,1}\ldots, \ldots
q_{\ell-m,\ell-m+n}$.  See the proof of Proposition 2.4.10.

(d) For $m=0,\ldots,\ell-1,n=0,\ldots,\ell-m-1$,
$$
[X^{\ell+m+n},\xi_{\ell-m}]_{c(m,n)} = 0, \eqno(2.24)
$$
with
$$
\eqalign{& c(m,n) = q_{\ell-m,1} \ldots, \ldots
q_{\ell-m,\ell-m-n},~~n\not= 0 \cr
&c(m,0) = q_{\ell-m,1} \ldots, \ldots q_{\ell-m,\ell-m+1}. \cr}
$$

\b \no {\bf Proof of (d).} The special case $m=n=0$ is
included in Part (a), and we check that
$c(0,0)=a(\ell,\ell)$.  We shall verify that all the
derivatives of (2.24) vanish.  First, when $n=0$, there is
only one nontrivial case, namely
$$
\eqalign{
\partial_{\ell-m}[X^{\ell+m},\xi_{\ell-m}]_{c(m,0)}
= &\,[\partial_{\ell-m}X^{\ell+m},\xi_{\ell-m}]_
{c(m,0)q_{\ell-m,\ell-m}}\cr
& +\bigl(q_{\ell-m,1} \ldots, \ldots q_{\ell-m,\ell-m}
-c(m,0)\bigr) X^{\ell+m}. \cr}\eqno(2.25)
$$
To show that this vanishes we calculate
$$\eqalign{
&\partial_{\ell-m}\partial_{\ell-m}
[X^{\ell+m},\xi_{\ell-m}]_{c(m,0)}\cr &=
(1+q_{\ell-m,\ell-m})\bigl(q_{\ell-1,1} ...,...
q_{\ell-m,\ell-m+1}-c(m,0)\bigr) X^{\ell+m}
= 0, \cr}
$$
$$
\eqalign{
&\partial_{\ell-m+1}\partial_{\ell-m}
[X^{\ell+m},\xi_{\ell-m}]_{c(m,0)}\cr
&\propto
\partial_{l-m+1}[X^{l+m-1},\xi_{l-m}]_{c(m,0)q_{l-m,l-m}}
\propto [X^{l+m-2},\xi_{l-m}]_{b(m,2)} = 0.
\cr}
$$
The statement is thus true for $n = 0$.  For $n\not= 0$
there are two items,
$$
\eqalign{\partial_{\ell-m}[X^{\ell+m+n},\xi_{\ell-m}]_{c(m,n)}
&= \bigl((q_{\ell-m,1}\ldots,\ldots q_{\ell-m,\ell-m-n}-c(m,n)\bigr)
X^{\ell+m+n} = 0, \cr
\partial_{\ell-m-n}[X^{\ell+m+n},\xi_{\ell-m}]_{c(m,n)} &\propto
[X^{\ell+m+n-1},\xi_{\ell-m}]_{c(m,n)q_{\ell-m-n,\ell-m}}. \cr}
$$
We note that $c(m,n)q_{\ell-m-n,\ell-m}=c(m,n-1)$.  When
$n\not= 1$ this depends on the constraint
$\sigma_{\ell-m,\ell-m-n}=1$.  When $n=1$, it makes use of
$\sigma_{\ell-m,\ell-m-1}q_{\ell-m,\ell-m}=1$, always valid.
The validity of (2.24) follows by induction on $n$.

\b \no {\it Corollary 2.4.12.} The element $X^{2\ell-1}$ is
a highest root vector for $C_\ell$.

\b
\no{\steptwo 3. Cohomology.}
\b

To each point $q$ in parameter space there corresponds a
free differential algebra $\cb_q$, an ideal ${\cal I}_q$
generated by the irreducible constants in $\cb_q$, and a
quotient algebra $\cb'_q = \cb_q/{\cal I}_q$.  In this
section we shall introduce a differential complex generated
by the constants.  The first three subsections are
tentative, exploratory and motivational.  The reader may
prefer to skip them.

\b \no{\bf 3.1.  The q-differential complex of a quantum
plane.}

\b Until further notice, a one-form is a map from
$\cb_q(1)$ to $\cb_q$ or to $\cb'_q$.

\b \no{\it 3.1.1.} Recall the notation $\sigma_{ij} =
q_{ij}q_{ji}, i,j = 1,...,N$.  Suppose first that
$\sigma_{ij} = 1, i\neq j = 1,...,N$, and that these are the
only constraints.  Then the polynomials
$$
C_{ij} = \sum_{k,l}C_{ij}^{kl}\xi_k\xi_l =
\xi_i\xi_j - q_{ji}\xi_j\xi_i,\quad i\neq j
= 1,...,N,
$$
are constants and the operators $D_q\hat C_{ij}$ are
identically zero,
$$
\hat C_{ij} :=\hat \partial_i\hat \partial_j -
q_{ji}\sl_j\sl_i,~~~ D_q\hat C_{ij} =\partial_i\partial_j -
q_{ji}\partial_j\partial_i  = 0,~~~i\neq j = 1,...,N.
$$
A one-form $y = (y_1,...,y_N)$ on $\cb_q(1)$ is called {\it
exact} if there is $x \in \cb_q$ such that $y_i = \partial_i
x, i = 1,...,N$, and it is said to be {\it closed} if
$$
C_{ij}(y):= \sum_{k,l}C_{ij}^{kl}\partial_ky_l =
\partial_iy_j -
q_{ji}\partial_jy_i = 0,~~~i\neq j = 1,...,N.
$$
It is usual to look upon the collection $\{C_{ij}(y)\}$ as
the components of an exact 2-form, but that is a point of
view that does not have a natural generalization.  Instead,
we shall say that a two-form $z = (z_{ij}), z_{ij} =
z(\xi_i,\xi_j)$, is {\it exact} if there is a 1-form $y$
such that $\sum_{k,l}C_{ij}^{kl}z_{kl} = C_{ij}(y), {i\neq j
= 1,...,N}$, and that it is closed if for all complex
coefficients $C^{ijk}$, such that
$\sum_{i,j}C^{ijk}\xi_i\xi_j =0=
\sum_{j,k}C^{ijk}\xi_j\xi_k$ in $\cb_q'$,
$\sum_{i,j,k}C^{ijk}\partial_i z_{jk} = 0$.  It is fairly
clear how this development can be completed to a
q-differential complex on quantum planes.  \b

\no{\it 3.1.2.} We consider the case when there is just one
condition on the parameters, $\sigma_{12} = 1$, and just one
constant $C_{12} = \sum_{k,l}C_{12}^{kl}\xi_k\xi_l = \xi_1\xi_2 -
q_{21}\xi_2\xi_1$.  Then
it makes sense to say that $y$ is closed if  $ C_{12}(y) =
\partial_1y_2 - q_{21}\partial_2y_1 =
0$, with no conditions on the other components, and that $z$
is exact if $\sum_{k,l}C_{12}^{kl}z_{kl} = C_{12}(y)$.  One
feels that, in this case as well, an associated differential
complex is lurking.

\bb \no {\bf 3.2.  The q-differential complex of a quantized
Kac-Moody algebra.}

\b \no{\it 3.2.1} The Serre constraints that define the
quantization of $A_2$ are
$$
\sigma_{12}q_{11} = \sigma_{12}q_{22} =1.
$$
The ideal ${\cal I}_q$ is generated by two constants in
$\cb_q(3)$,
$$
C_{112} = {1\over q_{12}}\xi_2\xi_1^2 - (1 + q_{11})\xi_1\xi_2\xi_1 + {1\over
q_{21}}\xi_1^2\xi_2,
$$
and $C_{221}$ defined analogously. Define
$$
\hat C_{112} :={1\over q_{12}}\sl_2\sl_1^2 -
(1 + q_{11})\sl_1\sl_2\sl_1 +
{1\over q_{21}}\sl_1^2\sl_2.,
$$
then the differential operator $D_q\hat C_{112}$ vanishes
identically.  A one-form $y = (y_1,...,y_N)$ may be called
exact if $y_i = \partial_i x, i = 1,...,N$ and it is natural
to say that it is strongly closed if it satisfies
$$
C_{112}(y) := \sum_{i,j,k}C_{112}^{ijk}
\partial_i\partial_jy_k = {1\over
q_{12}}\partial_2\partial_1y_1  -
(1 + q_{11})\partial_1\partial_2y_1 +
{1\over q_{21}}\partial_1^2y_2 = 0,
$$
as well as $C_{221}(y) = 0$.

\b \no{\bf 3.3.  q-differential complex in general.}

\b \no{\it 3.3.1.} If $C \in \cb_q$ is a constant, then
$\hat C := \widehat J C \in {\rm Ker}D_q$ and $C(y)$ is
obtained from the latter by replacing the right-most
$\hat\partial_i $ by $y_i$ and the other $\hat \partial_i$'s
by $\partial_i$'s operating on $y_i$.

\b \no{\bf 3.3.2.  Definition.  Proposition.} A one-form $y$
is said to be {\it exact} if there is $x \in \cb_q$ such
that $y_i = \partial_i x, i = 1,...,N$, and it is called
{\it closed} ({\it strongly closed}) if, $\forall C\in
\cb_q(2)$ ($\forall C\in\cb_q$), constant, $C(y) = 0$.
Every exact one-form is strongly closed, and therefore
closed, and every strongly closed one-form is exact.

\b \no{\it Proof.} Only the last statement needs
justification.  It is a corollary of Lemma 1.6.3, part (c).
For consider the expression
$$
x := \sum_{n \geq
1}\sum_{{\underline i}}A^{\underline i}\,\partial_{i_1}...
\partial_{i_{n-1}}y_{i_n}.
$$
Taking the derivatives of both sides of this equation we get
$$
\partial _kx - y_k =
\sum_{{\underline i}}\bigl(
\partial_kA^{ i_1...i_n} + \delta_k^{i_1}q_{i_1i_2}...q_{i_1i_n}
A^{i_2...i_n} \bigr)\partial_{i_1}...\partial_{i_{n-1}}y_{i_n}.
$$
This is equal to zero if $y$ is strongly closed, as can be
checked by consulting subsection 1.6.3.

\b \no{\it 3.3.3.  Definition.} A homogeneous constant (a
constant that is a sum of reorderings of a monomial,
Definition 1.2.1) is {\it irreducible} if it does not belong
to the ideal generated by constants of lower order.  A {\it
reducible} constant is a sum of polynomials of total degree
$k$ each of which contains a constant factor of lower
degree.

\b It is clear that the ideal ${\cal I}_q$
generated by all the constants is generated by the
irreducible constants.  It is easy to check that a one-form
$y$ is strongly closed iff $\forall C\in \cb_q$, $C$ an
irreducible constant, $C(y) = 0$.

\b \no{\it 3.3.4.} Before proceeding to complete the
construction of a differential complex in the general case
we present an example of a different kind, beyond the
context of Kac-Moody algebras.  Suppose at first that the
parameters satisfy the condition $\sigma_{123} :=
\sigma_{12}\sigma_{13}\sigma_{23} = 1$, but are otherwise in
general position.  Then there is just one irreducible
constant, and it can be written as follows
$$
C_{123} = {1\over q_{12}}\biggl(\xi_2(\xi_3\xi_1 +
\sigma_{12}q_{13}\xi_1\xi_3) -
q_{32}q_{12}(\xi_3\xi_1 +
\sigma_{12}q_{13}\xi_1\xi_3)\xi_2\biggr) + {\rm
cyclic}.
$$
Now suppose that
$$
\sigma_{12}\sigma_{23}\sigma_{13} =: \sigma_{123} =
\sigma_{124} =
\sigma_{134} = \sigma_{234} = 1.
$$
This implies that either $\sigma_{ij} = \sigma _{kl}$ or
else $\sigma_{ij} = -\sigma _{kl}$, for $i,j,k,l$ all
different.  Now there are 4 irreducible constants.  One can
verify that, if $\sigma_{ij} = \sigma _{kl}$, then there is
an identity
$$
\eqalign{
&\hat\partial_4\hat C_{123} + q_{14}q_{13}q_{34}
\hat\partial_1\hat C_{234}
+ {q_{21}q_{13}\over
q_{42}q_{43}}\hat\partial_2\hat C_{134}
+q_{34}q_{32}q_{24}\hat\partial_3\hat C_{124}\cr
&-{q_{24}\over q_{31}}\hat C_{412}\hat\partial_3 -
{q_{34}\over q_{12}}\hat
C_{423}\hat\partial_1 -
q_{13}q_{32}q_{34}\hat C_{413}\hat\partial_2 -
q_{14}q_{24}q_{34}\hat
C_{123}\hat\partial_4 = 0.\cr}
$$
Consequently, for any one-form $y$,
$$
\partial_4 C_{123}(y) + q_{14}q_{13}q_{34} \partial_1
C_{234}(y) +
{q_{21}q_{13}\over
q_{42}q_{43}} \partial_2  C_{134}(y) +q_{34}q_{32}q_{24}
\partial_3
C_{124}(y) = 0,
$$
or, if $(dy)_i = C_{jkl}(y)$ for each cyclic permutation
$i,j,k,l$ of $\{1,2,3,4\}$,
$$
\partial_4 (dy)_4  + q_{14}q_{13}q_{34} \partial_1 (dy)_1
+ {q_{21}q_{13}\over
q_{42}q_{43}} \partial_2  (dy)_2  +q_{34}q_{32}q_{24}
\partial_3  (dy)_3  = 0,
$$
The equations $dy_i = z_i, i = 1,2,3,4$ are not integrable
unless the two-form $z$ satisfies
$$
\partial_4 z_4  + q_{14}q_{13}q_{34} \partial_1 z_1  +
{q_{21}q_{13}\over
q_{42}q_{43}} \partial_2 z_2  +q_{34}q_{32}q_{24}
\partial_3 z_3  = 0;
$$
we may ask whether this condition is sufficient to guarantee
that $z$ can be expressed as $dy$.  Complete answers to all
these problems of integrability will now be found within a
study of the Hochschild cohomology of $\cb'_q$.

\b
\no{\bf 3.4. Hochschild complex.}
\b
Here we shall see that the idea of constructing a
differential complex based on $\cb'_q$ can be realized
within the setting of the ordinary Hochschild complex of
$\cb'_q$.

\b \no{\it 3.4.1.} Let $\ca\,$ be an associative {\bf
C}-algebra.  On $p$-chains $a_1\otimes ...  \otimes a_p \in
C_p(\ca), a_i \in {\cal E}, p\geq 1$, define a linear
boundary operator $\partial:C_p \rightarrow C_{p-1}, p \geq
2$, and $C_1 \rightarrow 0$, by
$$
\partial(a_1\otimes ... \otimes a_p) =
\sum_{i=1}^{p-1}(-)^{i+1}(a_1\otimes
a_2\otimes ...\otimes a_{i-1}\otimes
a_ia_{i+1}\otimes ...
\otimes a_p).
$$
One has $\partial\circ \partial = 0$.
\b
\no{\it 3.4.2. Examples.}
$$\eqalign{&\partial(a_1) = 0,\cr
&\partial(a_1\otimes a_2) = a_1a_2,\cr
&\partial(a_1\otimes a_2\otimes a_3) =
(a_1a_2\otimes a_3) - (a_1\otimes
a_2a_3).\cr}
$$
For $a \in C_p(\ca)$ we use the notation
$$
a = \sum a_{(1)}\otimes a_{(2)} \otimes a_{(3)}
\otimes...\otimes a_{(p)} =
\sum a_{(1)} \otimes a_{(2,...,p)}.
$$
The formula in 3.4.1 applies with indices in parentheses and
summation, for example $\partial(\sum a_{(1)}\otimes
a_{(2)}) =\sum a_{(1)}a_{(2)}$.  The summation sign will
usually be omitted.

\b \no{\bf 3.4.3.  Proposition.} The ${\bf C}$-algebra
$\cb(+) = \sum_{n \geq 1}\cb(n)$, freely generated by
$\xi_1,...,\xi_N$, has homology quotients $H_1(\cb(+)) \sim
\cb(1)$ and $H_p(\cb(+)) = 0, p \geq 2$.

\b \no{\it Proof.} All 1-chains are closed, and all
homogeneous 1-chains of degree higher than 1 are of the form
$a = \sum_i\xi_ia^i = \partial(\sum_i\xi_i \otimes a^i)$,
so
$H_1(\cb(+)) \sim \cb(1)$.  Every $a \in C_{p+1}(\cb(+)),~ p
\geq 1$,
$$
a = a_{(0)} \otimes a_{(1)} \otimes ... \otimes a_{(p)},
~~~ a_{(i)} \in
\cb(+),~~~ i = 0,1,...,p,
$$
is homologous to a chain of the form
$$
\hat a =\sum_{i=1}^N \xi_i \otimes b^i  \in C_{p+1}(\cb(+)),
~~~b^i =b^i_{(1)}\otimes ... \otimes b^i_{(p)}\in C_p(\cb(+)).
$$
Hence $\partial a = \partial\hat a =
\partial\bigl(\sum_i\xi_i\otimes b^i\bigr) =
\sum_i(\xi_ib^i_{(1)})\otimes b^i_{(2)}\otimes ...\otimes
b^i_{(p)} -\sum_i\xi_i\otimes \partial b^i$.  The degree of
$\xi_ib^i_{(1)}$ is greater than 1; therefore, if
$\partial a = 0$, then $\xi_ib^i_{(1)} = 0$.  Since there
are no relations in \cb, $b^i = 0$ and $a$ is homologous to
zero.

\b
As a corollary of the proof we have

\b \no{\bf 3.4.4.  Proposition.} One has $H_1(\cb'_q(+)) \sim
\cb(1)$ and $H_p(\cb'_q(+)) \sim \{\sum_i\xi_i \otimes b^i
\in Z_p(\cb'_q)\}, p\geq 2$.

\b\no{\it 3.4.5.} We introduce the spaces $C^p$ of linear
maps from $C_p$ to a left \ca-module $M$ and set $C^0 = M$.
The action of $a \in \ca$~ in $M$ will be denoted $\pi(a)$.
A linear coboundary operator $d: C^{p-1} \rightarrow C^p$ is
defined by
$$
d\tau(a)  = \pi(a_1)\tau(a_2\otimes...\otimes a_p) -
\tau(\partial a),~~~ a
=  a_1\otimes ... \otimes a_p.
$$
One has $d\circ d = 0$.

\b\no{\it 3.4.6. Examples.}
$$\eqalign{
&d\tau(a_1) = \pi(a_1)\tau,\cr
&d\tau(a_1\otimes a_2) = \pi(a_1)\tau(a_2) - \tau(a_1a_2),\cr
&d\tau(a_1\otimes a_2 \otimes a_3) =\pi(a_1)\tau(a_2\otimes a_3) -
\tau(a_1a_2 \otimes a_3) + \tau(a_1\otimes
a_2a_3).\cr}
$$
The linear extension takes the form
$$\eqalign{
&a \in C^1,~~d\tau(a) = \pi(a)\tau,\cr
&a \in C^2,~~d\tau(a) = \sum \pi(a_{(1)})\tau(a_{(2)}) - \tau(\partial a),\cr
&a \in C^3,~~ d\tau(a) = \sum \pi(a_{(1)}) \tau(a_{(23)})- \tau(\partial
a).\cr}
$$

The homomorphism $J_q\hskip-1mm: \cb_q \rightarrow \cb_q^*$,
generated by $\xi_i \mapsto \pi(\xi_i) = \partial_i\,$,
provides a new action of \cb$_q\,$ on $\cb_q\,$ and on
$\cb_q'$.  From now on we take $M = \cb_q \,$ or $M =
\cb_q'\,$ and write $\pi(a)$ for this action.  In the next
proposition $M = \cb_q$.

\b \no{\bf 3.4.7.  Proposition.} The Hochschild cohomology
of $\cb_q(+) := \bigoplus_{n\geq 1} \cb_q(n)$, with values
in $\cb_q$, vanishes except that $H^0(\cb_q(+),\cb_q)$ is
the linear span of $\cb(0)$ with the space of constants in
$\cb_q$.

\b \no{\it Proof.} For $p \geq 1$, let $z $ be a closed
$p$-cochain; we try to find a $(p-1)$-cochain $y$ such that
$$
z(a) = dy(a) = \pi(a_{(1)})y(a_{(2...p)}) - y(\partial a).
$$
When $p = 1$ the existence of $y$ is assured by Proposition
3.3.2, so from now on suppose that $p >1$.  We interpret the
equation, recursively in the degree, as the definition of
the last term, where the argument $\partial a$ has the
highest degree.  The obstruction is $\partial a = 0, a \neq
0$.  If all $a_{(i)}, i = 1,...,p$ are of degree
1 then, since there are no relations, $\partial a = 0$
implies that $a = 0$; this establishes the basis for the
recursion.  In general, when $\partial a = 0$, by
Proposition 3.4.3, there is a $(p+1)$-chain $b$ such that $a
= \partial b$ and we need to satisfy
$$
z(\partial b) = \pi(b_{(0)}b_{(1)})y(b_{(2...p)}) - \pi(b_{(0)})y(\partial
b_{(1...p)}),
$$
Since $z$ is closed, $z(\partial b) =
\pi(b_{(0)})z(b_{(1...p)})$, this holds if
$$
 z(b_{(1...p)}) = \pi(b_{(1)})y(b_{(2...p)}) -
 y(\partial b_{(1...p)}).
$$
This throws the solution of the equation $z(a) = dy(a)$ back
on the solution of the same equation with $a$ replaced by
$p$-chains of lower degree.  The existence of a base for the
recursion was established, and this completes the proof of
the theorem.

\b \no{\bf 3.4.8.  Theorem.} The Hochschild cohomology of
$\cb'_q(+) := \bigoplus_{n \geq 1}\cb'_q$, with values in
$\cb_q$ or in $\cb_q'$, vanishes except that
$H^0(\cb_q'(+),\cb_q)$ is the union of $\cb(0)$ and the
space of constants and $H^0(\cb_q'(+),\cb_q') = \cb(0) =
{\bf C}$.

\b \no{\it Proof.} We begin as in the proof of 3.4.7, but
$\partial a = 0$ no longer implies that $a$ is exact.  By
Proposition 3.4.4 every closed $p$-chain is homologous to
one of the form $a = \xi_i \otimes b^i$, so we have to show
that, when $z$ is closed, there is $y$ such that
$$
z(\xi_i \otimes b) = \partial_iy(b) -
y(\partial (\xi_i \otimes b)).\eqno(3.1)
$$
We need a lemma.

\b \no{\it 3.4.9.  Lemma.} Let $\{x_\alpha\}$ be any finite
collection of homogeneous elements of $\cb'_q$, all of the
same degree in each variable, satisfying a linear functional
relation $\sum_{\alpha} A^\alpha (x_\alpha) = 0$.  Then
there is a family $\{\Pi^\alpha\}$ of differential
operators, such that $\sum_{\alpha}A^\alpha(x_\alpha) =
\sum_{\alpha}\Pi^\alpha x_\alpha$.

\b \no{\it Proof.} This is just the statement that the
algebra $\cb^*_q$ of differential operators is the algebraic
dual of $\cb'_q$, see Theorem 1.4.5.

\b \no{\it 3.4.10.  Proof of the Theorem.} To settle the
integrability of Eq.(3.1), note that it is a linear relation
in a finite dimensional vector space.  Given the left side,
(3.1) has a solution $y$ if and only if the left side
satisfies all linear functional relations that hold for the
right hand side identically in $y$.  By the lemma, such
relations take the form
$$
\sum_{\alpha,i} \pi(c^{i,\alpha})\bigl[\partial_iy(b_\alpha) -
y(\partial(\xi_i \otimes b_\alpha))\bigr] = 0,
$$
where $\{b_\alpha\}$ is a family of $(p-1)$-chains and
$\{c^{i,\alpha}\}$ is a family of elements of $\cb'_q$.
Since this is required to hold identically for all $y$, both
terms must vanish separately,
$$
\sum_{\alpha,i} \pi(c^{i,\alpha}\xi_i)y(b_\alpha) = 0,~~
\sum_{\alpha,i}\pi(c^{i,\alpha})y(\partial (\xi_i \otimes
b_\alpha))\bigr] = 0.
$$
The second relation is equivalent to
$$\sum_{\alpha,i}c^{i,\alpha} \otimes
\partial(\xi_i \otimes b_\alpha) = 0;
$$
the first relation is satisfied for all $y$ if and only if
the differential operators are 0; hence (3.1) is integrable
if and only if
$$
\sum_{\alpha,i} \pi(c^{i,\alpha})z(\xi_i \otimes b_\alpha) = 0,
$$
for all families $\{b_\alpha\}$ and all $\{c^{i,\alpha}\}$,
such that $\sum_{\alpha,i}c^{i,\alpha} \otimes
\partial(\xi_i \otimes b_\alpha) = 0$ and $\sum_i
c^{i,\alpha}\xi_i = 0$; in other words for all $\{b_\alpha,
c^{i,\alpha}\}$ such that
$\sum_{\alpha,i} c^{i,\alpha} \otimes \xi_i \otimes
b_\alpha$ is closed.  But that is true if $z$ is closed.
The theorem is proved.

\b \no {\bf 3.5 Serre cohomology.}

\b \no{\it 3.5.1.} A zero-cochain $x \in C^0(\cb'_q,\cb'_q)$
is an element $x$ of ${\cal{B}}_q^\prime$, it is exact only
if $x=0$, and it is closed if $\pi(a)x=0,~\forall a\in
{\cal{B}}_q^\prime$, which (because there are no constants
in $\cb'_q$) is true iff $x\in {\cal{B}}(0)={\bf C}$.
\vskip.5cm \b \no{3.5.2.} A one-cochain $y$ on
${\cal{B}}_q^\prime$ is exact if there is
$x\in{\cal{B}}_q^\prime$ such that $\forall a\in
{\cal{B}}_q^\prime$,
$$
y(a)=\pi(a)x~. \eqno(3.2)
$$
\no Let $y_1$ denote the restriction of $y$ to $\cb_q(1)$,
then if $y$ is exact we have
$$
y_1(\xi_i) = \partial_ix~, \eqno(3.3)
$$
in which case we say that $y_1$ is exact.  Conversely, if
$y_1$ is an exact one-form on $\cb_q(1)$, then there is a
unique, exact one-cochain $y$ on ${\cal{B}}_q^\prime$ that
restricts to $y_1$ on $\cb_q(1)$. \footnote*{There are no constants in
$\cb_q(1)$; hence $\cb'_q(1) = \cb_q(1)$.}

\b \no {\it 3.5.3.} A one-cochain $y$ on
${\cal{B}}_q^\prime$ is closed if $\forall a,b\in
{\cal{B}}_q^\prime$,
$$
dy(a\otimes b) = \pi(a)y(b)-y(ab) = 0~. \eqno(3.4)
$$
This implies that, if $c=\sum_i a_i\xi_i$, then
$$
\sum_i \pi(a_i)y(\xi_i) = y(c)~, \eqno(3.5)
$$
and in particular that
$$
c=\sum_i a_i\xi_i=0 \Rightarrow c(y_1):=
\sum_i \pi(a_i)y_1(\xi_i)=0~. \eqno(3.6)
$$
Compare 3.1.1.  Conversely, if $y_1$ satisfies (3.6), then a
unique, closed one-cochain $y$ on ${\cal{B}}_q^\prime$ is
determined by solving (3.4), recursively in the degree of
$ab$.

\b \no {\bf Definition 3.5.4}.  A $p$-form $z_1$ on
$\cb_q(1)$, with values in ${\cal{B}}_q^\prime$, is exact if
there is a $(p-1)$ -form $y_1$ on $\cb_q(1)$ such that
$$
z_1(a) = \pi(a_{(1)})y_1(a_{(2)}\otimes\ldots\otimes a_{(p)})~,
\quad
\forall a=a_{(1)}\otimes\ldots\otimes a_{(p)}~, ~~
\partial a=0~. \eqno(3.7)
$$
\no It is closed (strongly closed) if
$$
\pi(a_{(0)})z_1(a_{(1)}\otimes\ldots\otimes a_{(p)})=0~,
\quad
\forall  a = a_{(0)}\otimes a_{(1)}\otimes\ldots\otimes
a_{(p)}~, \quad
\partial a = 0~, \eqno(3.8)
$$
where $a_{(1)},\ldots,a_{(p)}\in \cb_q(1)$ in both
cases, and $a_{(0)}\in {\cal{B}}_q(1)$ ($a_{(0)}\in
{\cal{B}}_q^\prime$).

\b \no {\it 3.5.5} It is obvious that, if $z$ is an exact
$p$-cochain on ${\cal{B}}_q^\prime$, then its restriction
$z_1$ is an exact $p$-form on $\cb_q(1)$.  Conversely, if
$z_1$ is an exact $p$-form on $\cb_q(1)$, expressed in terms
of a $(p-1)$-form $y_1$ as in (3.7), let $y$ be any
$(p-1)$-cochain on ${\cal{B}}_q^\prime$ that extends $y_1$;
then the formula $z=dy$ extends $z_1$ to an exact
$p$-cochain on ${\cal{B}}_q^\prime$.

\b \no {\it 3.5.6.} If $z$ is a closed $p$-cochain on
${\cal{B}}_q^\prime$ then its restriction $z_1$ is a
strongly closed $p$-form on $\cb_q(1)$.  Conversely, if
$z_1$ is a strongly closed $p$-form on $\cb_q(1)$, then it
extends to a closed $p$-cochain on ${\cal{B}}_q^\prime$.  To
show this consider the condition
$$
dz(a):= \pi(a_{(0)})z(a_{(1)}\otimes\ldots\otimes a_{(p)})
-z(\partial a)=0
\eqno(3.9)
$$
for $a=a_{(0)}\otimes a_{(1)}\otimes\ldots\otimes a_{(p)}$,
with $\hbox{deg}(a_{(i)}) \geq 1, i = 1,...,p$ and of total
polynomial degree
$$
n(a)=\sum^p_{i=0} \hbox{deg}(a_{(i)})\geq p+1~.
$$
When $n(a)=p+1$, this condition determines $z(\partial a)$
in terms of $z_1$; this amounts to a partial determination
of $z(b)$ for $n(b)=p+2$.  Proceed recursively in the total
degree $n(a)$.  Eq.(3.9) is solvable for $n(a)=p+1$.
Suppose it can be solved for $n(a)=p+1,p+2,\ldots,p+k$, and
let $n(a)=m=p+k+1$.  Then $\partial a$ is of order $m+1$,
and (3.9) can be solved for $z(\partial a)$ provided
$$
\pi(a_{(0)})z(a_{(1)}\otimes\ldots\otimes a_{(p)})=0~, \quad
\forall a = a_{(0)}\otimes ...
\otimes a_{(p)},~\partial a=0~. \eqno(3.10)
$$
The apparent obstruction to (3.10) comes from the fact that 
$z(...)$ on the left hand side is already known for the case 
that the argument is exact.  Setting this argument equal to 
$\partial b$, and using (3.9) for lower degrees, we reduce 
the left hand side of (3.10) to 
$\pi(a_{(0)}b_{(0)})z(b_{(1...p)})$.  But this is the same 
as $z(\partial a)$ and hence zero by induction in the degree 
of $a_{(1...p)}$.

We conclude as follows.

\b \no {\bf Theorem 3.5.7.} Let ${\cal D}$ denote the
differential complex of which the cochains are
$\cb_q(1)$-forms restricted to closed chains, and the
differential is given, for cochains valued in $\cb'_q$ by
the formula
$$
dz(a_{(0)}\otimes a_{(1)}\otimes\ldots\otimes a_{(p)}) =
\pi(a_{(0)})z(a_{(1)}\otimes\ldots\otimes a_{(p)})~. 
\eqno(3.11)
$$
The cohomology is nontrivial if $\cb_q$ admits irreducible
constants\footnote*{See Definition 3.3.3.} of order higher
than 2.  If
$$
C =\sum_{\underline i}C^{\underline i}\,
\xi_{i_1}...\xi_{i_k}
$$
is a constant, then
$$
\sum_{\underline i}C^{\underline i}\,
\partial_{i_1}...\partial_{i_k}  = 0,
$$
and $z$ is exact iff $z$ is strongly closed; that is, iff
for all constants $C$,
$$
\sum_{\underline i}C^{\underline i}\,\partial_{i_1}...
\partial_{i_{k-1}}z(\xi_{i_k}\otimes\xi_{j_1}\otimes ...\otimes
\xi_{j_p}) = 0,
$$
(If $x$ is a closed 0-form, then $\partial_ix=0$, so $x\in
{\cal{B}}_0 = {\bf C}$, and no 0-cochain is exact, so
$H^0={\bf C}$.)

\b The restriction to closed chains is very
natural; in the ($q$-) commutative case it restricts the
cochains to be ($q$-) alternating.  But if there are no
constants (in $\cb_q$) of order two, then the restriction is
moot.

\bb \no 
{\steptwo{4.  Classification of Constraints.}} 
\b
\no{\bf 4.1.  The Constants in ${\cal B}_{1...n}$.} 
\b \no
{\it 4.1.1.} From now on we drop the suffix $q$ on $\cb$ and
denote by ${\cal B}_n$ the space of polynomials of order
$n$.  For any subset $s = \{i_1,...,i_{|s|}\}\subset
\{1,...,n\}$, with cardinality $|s|$, denote by $\cb_{(s)}$
the subspace of $\cb_{|s|}$ consisting of all polynomials of
degree $|s|$ that are separately linear in
$\xi_{i_1},...,\xi_{i_{|s|}}$, and let
$$
\sigma_{(s)} = \prod_{i\neq j \in s}q_{ij}.\eqno(4.1)
$$

\b \no {\bf 4.1.2.  Theorem.} Assume that there are no
constants in $\cb_{(s)}$, for any proper subset $s \subset
\{1,...,m\}$.  Then the dimension of the space of constants
in ${\cal B}_{1...m}$ is
$$
\matrix{ (m-2)!&{\rm if }~\sigma_{1...m} = 1,\cr &\cr 0&{\rm
otherwise.}\cr}\eqno(4.2)
$$

\b \no {\it Proof.} By induction.  The statement is true for
$m = 2$.  We assume that it is true for $m = 2,...,n-1$ and
prove that it is true for $m=n$.  Since there are no
constants in ${\cal B}_{(s)}$ and the statement is true for
$m = 2,...,n-1$ it follows that $\sigma_{(s)} \neq 1$, for
any proper subset $s$ of $1...n$.

\b \no {\it 4.1.3.} Let $X =\sum X(1...n)\,\xi_1...\xi_n$
~(sum over all permutations of $12...n$) be a constant in
${\cal B}_{1...n}$.  The equations $\partial_iX = 0, i =
1,...,n$ are
$$
\eqalign{
&X(n12...) + q_{n1}X(1n2...) + q_{n1}q_{n2}X(12n3...) +
... \cr
&\phantom{X(n12...)} + q_{n1}...q_{nn''}X(12...n''nn') +
q_{n1}...q_{nn'}X(12...n'n) = 0,\cr
&X(1n2...) + q_{1n}X(n\underline{12...n'}) = 0,\cr
&X(\underline{12}n3...) +
q_{12}q_{1n}X(2n\underline{13...n'}) = 0\cr & ... \cr
&X(\underline{12...n'}n) + q_{12}...q_{1n}X(23...n'n1) = 0.
\cr}\eqno(4.3)
$$
and those obtained from this set by permuting the indices
$1,...,n'$.  We write $12...n'$ instead of $i_1...i_{n'}$
for any permutation of $1,...,n'$.

\b \no {\it 4.1.4.  Conventions.} We have used $n' = n-1,
n'' = n-2,$ and
$$
\eqalign{
X(...\underline{1...k}...) &:= S(1...k)X(...1...k...)\cr
&:= X(...1...k...) + q_{12}X(...21...k...) + ... +
q_{12}...q_{1k}X(...23...k1...),
\cr}
$$
 In addition, set
$k' = k+1, k'' = k+2$ and
$$
X(\overline{1...p}...mn...) = S(p...m)...S(2...m)S(1...m)X(1...p...mn...).
$$
Here $m$ is by definition the index that stands to the left
of $n$.  Note that the permutations act on the symbols, not
on the spaces, and that they affect the parameters, for
example $S(12)q_{13} = q_{23}S(12)$.  Finally, a lot of
typing is saved by introducing
$$
Z(1...mn...) := q_{n1}...q_{nm}X(1...mn...).
$$

\b \no {\it 4.1.5.  Remark.} The operator that sends
$X(...i_1...i_k...)$ to $X(...\underline{i_1...i_k}...)$ for
each permutation $i_1...i_k$ of $1...k$ corresponds to
differentiation; $X(...\underline{1...k}...)$ is the
coefficient of $\xi_2...\xi_k$ in
$\partial_1\sum_{\underline i}
X(...i_1...i_k...)\xi_{i_1}...\xi_{i_k}$.  It is invertible
if and only if $\partial_iX = 0, i = 1,...,k$ implies that
$X = 0$.

\b \no{\it 4.1.6.} With this notation equations (4.3) take
the form
$$
\eqalign{
&Z(n12...) + Z(1n2...) + Z(12n...) + ... + Z(12...n'n) = 0,\cr
&Z(\underline{1...k}n...n') +
q_{12}...q_{1k}\sigma_{1n}\,Z(2...kn\underline{1k'...n'})
= 0,
\quad k =1,...n'.
\cr}\eqno(4.4-5)
$$
Applying $S(1...n')$ to the long equations and using the
short ones to push the index $n$ towards the right we obtain
$$
(1-\sigma_{1n})\bigl(Z(1n2...n') + Z(\underline{12}n3...n')
+ ... +
Z(\underline{12...n'}n)\bigr)  =  0.
$$
This completes {\it Step 1}.  By stipulation $\sigma_{1n}
\neq 1$, so that the first factor can be dropped; hence
$$
\sum_{p=1}^{n'}Z(\underline {1...p}np'...n') = 0.
$$
Next apply $S(2...n')$ to this equation and use the short
equations (4.5) in the same way.  We claim that the result
after {\it Step} $k$ is
$$
\biggr(1-{1 \over q_{k1}...q_{kk-1}\sigma_{kn}}
P_{k...1}\biggl)\sum_{p=k}^{n'}
Z(\overline{1...k}k'...pnp'...n')
= 0,\eqno(4.6)
 $$
where $P_{k...1}$ is the cyclic permutation that takes $m
\mapsto m-1 $ $ (m = 2,...,k)$, $1 \mapsto k$.  To verify
this we carry out the next step.  We need a simple Lemma.

\b \no {\it 4.1.7.  Lemma.} The first factor in (4.6) is
invertible if and only if $ \sigma_{1...kn} \neq 1$.

\b \no {\it Proof of the Lemma.} The permutation $P_{k...1}$
is of order $k$.  Iteration of (4.6) leads to $(1-A)\sum Z =
0$, with
$$
A := \biggr({1 \over
q_{k1}...q_{kk-1}\sigma_{kn}}P_{k...1}\biggl)^k =
\sigma_{1...kn}.
$$
Hence, if $A \neq 1$, then (4.6) implies that $\sum Z = 0$.

\b \no{\it 4.1.8.} The first factor in (4.6) can therefore
be dropped.  Next, apply $S(k'...n')$ to get
$$
\sum_{p=k'}^{n'}Z(\overline{1...k'}k''...pnp'... n') +
\sum_{p=k}^{n''}q_{k'k''}...q_{k'p'}
Z(\overline{1...k}k''...p'n\underline{k'p''.
..n'})
=0.\eqno(4.7)
$$
The second term is
$$
\sum_{p = k'}^{n'}QZ(\overline{1...k}k''...pn
\underline{k'p'...n'})
= \sum_{p = k'}^{n'}\hat Q
P_{k'...1}Z(\overline{2...k'}k''...pn{1p'...n'}),
\quad \hat
Q = q_{k'k''}...q_{k'p}.
$$
Using the short equations we transform the summand to
$$\eqalign{
-QP_{k'...1}[q_{12}...q_{1p}&\sigma_{1n}]^{-1}
Z(\overline{1...k'}k''...pnp'...n'
)\cr
&= -[q_{k'1}...q_{k'k}\sigma_{k''n}]^{-1}
P_{k'...1}Z(\overline{1...k'}k''...pnp'...n').
\cr}
$$
Thus (4.7) becomes
$$
\biggr(1-{1
\over
q_{k'1}...q_{k'k}\sigma_{k'n}}P_{k'...1}\biggl)
\sum_{p=k'}^{n'}Z(\overline{1...k'}...pnp'...n') = 0.
\eqno(4.8)
$$
This proves our claim (4.6).  The process ends with $k =
n'$, when (4.6) reduces to
$$
\bigg(1 -
{1\over
q_{n'1}...q_{n'n''}\sigma_{n'n}}\,P_{n'...1}\bigg)
\,Z(\overline{1...n''}n'n) =0
.\eqno(4.9)
$$
By Lemma 4.1.3 the first factor is invertible unless
$\sigma_{1...n} = 1$.

\b \no{\it 4.1.9.} We have shown that, if $\sigma_{1...n}
\neq 1$, then $Z(\overline{1...n''}n'n)$ must vanish, along
with all the coefficients related to it by permutations of
the indices $1...n'$.  By Remark 4.1.5 it follows that
$Z(i_1...i_{n'}n) = 0$ for any permutation $i_1...i_{n'}$ of
$1...n'$.  In the same way we use the short equations (4.5)
to show that in that case all the coefficients are zero.
Thus, when $\sigma_{1...n} \neq 1$ there are no constants in
${\cal B}_{1..n}$.  Conversely, if $\sigma_{1...n} = 1$,
then Eq.s (4.9) has nontrivial solutions.  Now Eq.(4.9)
fixes the ratios of $n'$ coefficients related by
$P_{n'...1}$.  Among all the coefficients
$Z(i_n...i_{n'}n)$, $(n-2)!$ remain arbitrary.  Choosing any
solution of (4.9) we use the short equations to construct a
unique constant in ${\cal B}_{1...n}$.  The space of all
these constants has dimension $(n-2)!$ and Theorem 4.1.2 is
proved.

\b \no {\it 4.1.10.} As to determining the constants in
$\cb_{(s)}$ when the set $s$ contains repetitions, we do not
yet have complete results.

\b \no {\it 4.1.11.} We have nothing to say about the space
of constants in $\cb_{(s)}$ if $\sigma_{(r)} = 1$ for some
proper subset $r\subset s = \{1,...,n\}$.  That seems to
present a much more difficult problem.

\bb \no {\bf 4.2.  Constants of order 3.}

\b \no{\it 4.2.1.} We denote by ${\cal B}_{123}$ the
subspace of polynomials linear in $\xi_1,\xi_2,\xi_3$,
separately.  \b \no {\bf Proposition.} There are constants
in ${\cal B}_{123}$ iff
$$
(\sigma_{12}-1)(\sigma_{23}-1)(\sigma_{13}-1)(\sigma_{123} -1)
= 0.
$$

\b \no{\it Corollary.} Generically, there are no constants
and dim $\cb'_{(123)} = 6$.

\b \no{4.2.2.} Of special cases there are 5 essentially
different kinds.

\b (1) The constraint $\sigma_{123} = 1$. The space of constants is
1-dimensional with basis
$$
C_{123} = \bigl({1 \over q_{31}} - q_{13}\bigr)(\xi_1\xi_2\xi_3 +
q_{31}q_{32}q_{21}\xi_3\xi_2\xi_1) +  {\rm cycl.~ perm.}~.
$$
The intersection ${\cal I}_q \cap \cb_{123}$ is generated by
$C_{123}$ and the subspace $\cb'_{123}$ of the quotient is
5-dimensional.

\b (2) The constraint $\sigma_{12}=1$.  Then
$$
C_{12} = \xi_1\xi_2 -q_{21}\xi_2\xi_{1}
$$
is a constant.  The space of constants in $\cb_{123} $ is
1-dimensional with basis
$$
C_{12}\xi_3 - q_{31}q_{32}\xi_3C_{12}.
$$
The ideal ${\cal I}_q$ is generated by $C_{12}$, ${\cal
I}_q\cap \cb_{123}$ is two-dimensional and $\cb'_{123}$ is
4-dimensional.

\b (3) Two constraints: $\sigma_{12} = 1 =\sigma_{13}$.  The
space of constants in $\cb_{123}$ is two-dimensional and
$\cb'_{123}$ is two-dimensional with basis
$\{\xi_1\xi_2\xi_3, \xi_1\xi_3\xi_2\}$.

\b (4) Three constraints: $\sigma_{12} = \sigma_{13} =
\sigma_{23} = 1$.  The space of constants in $\cb_{123}$ is
the same as in the previous case, but $\cb'_{123}$ is
1-dimensional.

\b (5) The other case of two constraints, $\sigma_{12} = 1,
\sigma_{123} = 1$ yields a surprise.  When $\sigma_{12} =
1$, and $\sigma_{13} \neq 1 \neq \sigma_{23}$, then there
are no special cases: no further constant appears as
$\sigma_{123}$ takes the value 1.

\b \no {\it 4.2.3.  Summary.} We have the following complete
list of generators of ideals in $\cb_{123}$,
$$
\matrix{{\rm Constraints}&{\rm Generators}&
{\rm \#~ const.~}\in {\cal
B}_{123}&\dim {\cal I} \cap {\cal
B}_{123}\cr
\sigma_{12}\sigma_{23}\sigma_{13}=1& C_{123} &1&1\cr
\sigma_{12}=1& \xi_1\xi_2 - q_{21}\xi_2\xi_1 & 1&2\cr
\sigma_{12}=1,\sigma_{23}=1&\xi_1\xi_2 - q_{21}\xi_2\xi_1,
~\xi_2\xi_3 -
q_{32}\xi_3\xi_2
&2&4\cr
\sigma_{ij}=1,~ i \neq j &\xi_i\xi_j-q_{ji}\xi_j\xi_i\,, ~i \neq j&2&5
}
$$

\b \no {\it 4.2.4.} Constants of order 3, in two generators,
are of the Serre type.  Constants in one generator, $C =
\xi_1^3$ require $q_{11}^3 = 1$.  This completes the survey
of constants of order 3.

\bb
\no{\steptwo Appendix.  Alternative expressions for some
constants.}

\b The constants $C_{ij}, i\neq j$, given in 3.1.1 can be
written, up to a factor, as
$$
C_{ij}={\rm antisym}\bigl(\sqrt{q_{ij}}\;\xi_i\xi_j\bigr),
$$
where the antisymmetrizer is just $1-P_{ij}$, $P_{ij}$ being the
permutation that exchanges the indices $i,j$.

\b The constants $C_{ijk}, i\neq j\neq k$, given in 3.3.4
for the case 123, can be written, up to a factor, as
$$
C_{ijk}={\rm
sym}\bigl(\sqrt{q_{ij}q_{jk}q_{ki}}(q_{ki}^{-1}-q_{ik})
\xi_i\xi_j\xi_k\bigr),$$

where the symmetrizer is
$1+P_{ij}+P_{jk}+P_{ki}+P_{ij}P_{jk}+P_{ij}P_{ki}$.

These constants satisfy the identity given in 3.3.4, and
which we rewrite up to a factor as
$$
\bigl(\sqrt{q_{12}q_{13}q_{14}}\;\hat\partial_1\hat
C_{234}-\sqrt{q_{21}q_{31}q_{41}}\;
\hat C_{234}\hat\partial_1\bigr)+\,{\rm cycl.}=0.
$$

\bb \no {\steptwo Acknowledgements.} We thank Moshe Flato
for criticism.  A.G. thanks the Fundaci\'on Del Amo for
financial support and the Department of Physics of UCLA
for hospitality.

\ve
\no{\steptwo References.}

\item {[F]} C. Fr\o nsdal, Generalization and Deformations
of Quantum Groups, Publ.  Res.  Inst.  Math.  Sci.  {\bf
33}, 91-149 (1997) (q-alg/9606020).

\item {[FR]} I.B. Frenkel and N.Yu.  Reshetikhin, Quantum
Affine Algebras and Holonomic Difference Equations, Commun.
Math.  Phys.  {\bf 146} (1992) 1-60.

\item {[J]} M. Jimbo, A $q$-difference analogue of
$U(gl(N+1))$, Hecke algebra and the Yang-Baxter equation.
Lett.  Math.  Phys.  {\bf 11}, 247-252 (1986).

\item {[K]} M. Kashiwara,  On crystal basis of the Q-analogue of universal
enveloping algebras, Duke Math.  J. {\bf 63} (1991)
465--516.

\item {[L]} G. Lusztig, {\it Introduction to quantum groups}
Birkhaeuser, 1993.  Series: Progress in mathematics, vol.
110.

\item {[N]} H. Nikolai, A hyperbolic Kac-Moody algebra from
supergravity, Phys.Lett.  {\bf B276} (1992) 333--340.

\item {[S]} V. A. Smirnov, {\it Form Factors in Completely
Integrable Models of Quantum Field Theory}, Advanced Series
Math.  Phys., World Scientific 1992.

\item {[V]} A. Varchenko, {\it Multidimensional
hypergeometric functions and representation theory of Lie
algebras and quantum groups}, Advanced series in
mathematical physics; vol.  21 World Scientific 1995.  \ve

\end